\documentclass[final,leqno]{siamltex}
\title{On the extended randomized multiple row method for solving linear least-squares problems}
\author{Nian-Ci Wu\footnotemark[1]~,~\
        Chengzhi  Liu\footnotemark[2]~,~\
        Yatian Wang\footnotemark[3]~,
        \and
        Qian Zuo\footnotemark[4]\  \footnotemark[5]\  \footnotemark[6]}
\usepackage{amsmath}
\usepackage{amssymb}
\begin{document}
\maketitle
\renewcommand{\thefootnote}{\fnsymbol{footnote}}
\footnotetext[1]{School of Mathematics and Statistics, South-Central Minzu University, Wuhan 430074, China.}
\footnotetext[2]{School of Mathematics and Finance, Hunan University of Humanities, Science and Technology, Loudi 417000,  China.}
\footnotetext[3]{School of Mathematics and Statistics, Wuhan University, Wuhan 430072, China.}
\footnotetext[4]{Center on Frontiers of Computing Studies, Peking University, Beijing 10087, China.}
\footnotetext[5]{School of Computer Science, Peking University, Beijing 10087, China.}
\footnotetext[6]{Corresponding author. E-mail address: {{\sf zuoqian@pku.edu.cn} (Q.  Zuo).}}
\renewcommand{\thefootnote}{\arabic{footnote}}

\begin{abstract}
The randomized  row method is a popular representative of the iterative algorithm because of its efficiency in solving the overdetermined and consistent systems of linear equations. In this paper, we present an extended randomized multiple row method to solve a given  overdetermined and inconsistent linear system and analyze its computational complexities at each iteration. We prove that the proposed method can linearly converge in the mean square to the least-squares solution with a  minimum Euclidean norm. Several numerical studies are presented to corroborate our theoretical findings. The real-world  applications, such as image reconstruction and large noisy data fitting in  computer-aided geometric design, are also presented for illustration purposes.
\end{abstract}



\pagestyle{myheadings}
\thispagestyle{plain}
\markboth{The extended randomized multiple row method}{N.-C. Wu et al.}

\section{Introduction}
In a wide variety of scenarios, such as subroutines of tomographic imaging \cite{14AH} and data fitting in  computer-aided geometric design (CAGD) \cite{18LMD}, the research of efficient algorithms for solving large linear systems therein becomes an important topic. To be specific, we consider the linear system of the following form
\begin{align}\label{eq:Ax=b}
 A \bx = \bb,
\end{align}
where $A\in\Rc^{m\times n}$ is an $m$-by-$n$ ($m>n$) real coefficient matrix and $\bb\in\Rc^m$ is an $m$-dimensional right-hand side, and $\bx\in\Rc^n$ is an unknown $n$-dimensional vector to be solved.

With any $\bb\in {\Rc}^m$, we can uniquely write it as
$\bb= \bb_R + \bb_N$,  where $\bb_R:=AA^{\dag}\bb$ and $\bb_N:=(I_m-AA^{\dag})\bb$ are the projections of $\bb$ onto the range of $A$ (denoted by $\Rc(A)$) and the null space of $A^T$ (denoted by $\Nc(A^T)$), respectively, the symbol $\dag$ denotes the Moore-Penrose  pseudoinverse, and $I_m$ is the identity matrix of size $m$. In the computed tomography community, $\bb_R$ and  $\bb_N$ also name the vectors of exact data and corruptions, respectively. A challenge in the practical settings is that there is almost always corruption present in the large-scale real-world data,  or rather, $\bb_N\neq 0$. It may be caused by data collection, transmission, or modern storage systems. A common and widely studied approach is to seek the least-squares solution
$\bx_{LS} =\arg\min_{\bx\in\Rc^n} \BT{\bb - A\bx}$,
where $\|\cdot \|$ is used to represent the $2$-norm of either a vector or a matrix. It exists and is unique, and has the analytical expression $\bx_{LS} :=A^{\dag}\bb + \left(I_n - P_{A^T}\right)\by$ in which $\by$ is an $n$-dimensional vector and $P_{A^T}$ is the orthogonal projection onto $\Rc(A^T)$. By orthogonality it implies that $\BT{\bx_{LS}}=\BT{A^{\dag}\bb}+\BT{\left(I_n - P_{A^T}\right)\by}$ and follows that $\bx^{\ast} :=A^\dag\bb$ has the minimum  Euclidean norm \cite{96Bjo,13GL}. Many classic while effective iterative methods have been devised to compute $\bx^{\ast}$, e.g., the extended row algorithm \cite{95Pop,98Pop}.

Kaczmarz method is a well-studied row-oriented method due to its simplicity, efficiency, and intuitive geometric meaning \cite{Kac37}, which has been proved to be convergent to the least Euclidean norm solution of the consistent linear system. Popa utilized two Kaczmarz iterates and extended them to solve the inconsistent linear system \cite[Algorithm (R)]{95Pop}. At the $k$th iteration, the updates are given by
\begin{equation}\label{rek-zf}
\left \{
\begin{array}{l}
   \bx^{(k+1)}  =  \bx^{(k)} +  \frac{b_i - \by^{(k)}_i - A_{i,:}  \bx^{(k)}}{\BT{A_{i,:}}}A_{i,:}^T  \vspace{1.5ex}\\
   \by^{(k+1)}  =  \by^{(k)} - \frac{ A_{:,j}^T \by^{(k)}}{\BT{A_{:,j}}} A_{:,j}
\end{array}
\right.
\end{equation}
for $k=0,1,2,\cdots$, where $A_{i,:}$ and $A_{:,j}^T$ individually denote the $i$th and $j$th rows of $A$ and $A^T$, the indices $i\in [m]$ and $j\in[n]$ are chosen in the cycle orders, and the range of integers from $1$ to $\ell$ is written $[\ell]$. Zouzias and Freris separately selected
 the indices $i$ and $j$ with probabilities $\BT{A_{i,:}}/\BF{A}$ and $\BT{A_{:,j}}/\BF{A}$ and provided the randomized  extended Kaczmarz (REK) method  \cite{13ZF}.
Since $\by^{(k+1)}$ is a better approximation to $\bb_N$ \cite[Theorem 1]{19Du}, Du applied the Kaczmarz iterate to the linear system $A\bx = \bb - \by^{(k+1)}$ at the second-half step and obtained a slightly different REK iteration scheme.
In \cite{19BW}, Bai and Wu selected the index $j$ in a given cyclic order and proposed a partially randomized extended Kaczmarz method. Shortly afterwards, to further accelerate REK, they transformed the inconsistent system \eqref{eq:Ax=b} into a consistent augmented linear system, directly applies the greedy randomized Kaczmarz method \cite{18BW1} to it, and presented the greedy randomized augmented Kaczmarz method \cite{21BW}.

The selection of $A_{i,:}$ equals to use of an independent random variable $\bmu_i^T$ to multiply the left of $A$, where $\bmu_i$ denotes the $i$th standard unit coordinate vector with size $m$.  Wu and Xiang generalized it into matrix form and developed a framework for the analysis and design of REK. Some new iteration schemes were also given, such as the Gaussian extended Kaczmarz (GEK) method \cite{211WX}, defined by
\begin{equation}\label{gek}
\left \{
\begin{array}{l}
\by^{(k+1)}  = \by^{(k)} - \frac{\bze^T A^T \by^{(k)}}{\BT{A \bze}} A \bze \vspace{1.5ex}\\
\bx^{(k+1)}  = \bx^{(k)} + \frac{\bet^T(\bb - \by^{(k+1)} - A\bx^{(k)})}{\BT{A^T \bet}} A^T \bet,
\end{array}
\right.
\end{equation}
where the non-zero vectors $\bze$ and $\bet$ obey the standard normal distribution.

Variants of REK that exploit more than a single row at each iteration, often referred to as block methods, are extensively studied. Two particular methodologies have proven popular. The first one is the projective block method, in which each iteration is projected onto the subspace defined by multiple rows, such as the randomized double block Kaczmarz method \cite{15NZZ}. This kind of method is difficult to parallelize and needs to calculate the Moore-Penrose inverse. The second one is the pseudoinverse-free method, in which the projections of the previous iterate onto each row in a block are computed; see, e.g., the two-subspace randomized extended Kaczmarz method for a related fast solver \cite{22Wu} and the randomized extended average block Kaczmarz (REABK) for a survey on the latter \cite{20DSS}.

We call $\left\{\mathcal{I}_i\right\}_{i=1}^{s}$ a partition of $[m]$ if $\mathcal{I}_i \cap \mathcal{I}_j = \emptyset$ for $i\neq j$ and $\cup_{i=1}^s \mathcal{I}_i = [m]$. Similarly, let $\left\{\mathcal{J}_j\right\}_{j=1}^{t}$ be a partition of $[n]$.  REABK generates a new estimate $\bx^{(k+1)}$ via
\begin{equation}\label{reabk}
\left \{
\begin{array}{l}
\by^{(k+1)}  = \by^{(k)} - \alpha\frac{A_{:,\mathcal{J}_j}A_{:,\mathcal{J}_j}^T \by^{(k)}}{\BF{A_{:,\mathcal{J}_j}}} \vspace{1.5ex} \\
\bx^{(k+1)}  = \bx^{(k)} + \alpha\frac{A_{\mathcal{I}_i,:}^T \left(\bb_{\mathcal{I}_i}- \by^{(k)}_{\mathcal{I}_i} - A_{\mathcal{I}_i,:} \bx^{(k)} \right)}{\BF{A_{\mathcal{I}_i,:}}},
\end{array}
\right.
\end{equation}
where $i\in [s]$ (resp., $j\in [t]$) is picked with probability $\BF{A_{\mathcal{I}_i,:}}/\BF{A}$ (resp., $\BF{A_{:,\mathcal{J}_j}}/\BF{A}$), $A_{\mathcal{I}_i,:}$ (resp., $A_{:,\mathcal{J}_j}^T$) denotes the row submatrix of $A$ (resp., $A^T$) indexed by $\mathcal{I}_i$ (resp., $\mathcal{J}_j$), and $\alpha$ is a step size.

In this work, we give a block extended row method to deal with the general linear least-squares problems. Our approach is based on several ideas and tools, including an extended row iteration scheme and the promising random index selection strategy. The convergence analysis reveals that the resulting algorithm has a linear (sometimes referred to as {\sf exponential})  convergence rate, which is bounded by explicit expression. As with many iterative methods in the row family \cite{19BW, 19Du, 20DSS, 95Pop, 211WX, 22Wu, 13ZF}, our method relies on the information on the right-hand side. The underlying idea is that first using an additional row iterate to output an approximation  of $\bb_R$ and then utilizing another row iterate to an asymptotical  consistent  linear system. The theoretical results are derived in Section \ref{sec:ERMR}. Furthermore, we develop a short recursion to update the new approximation efficiently.
In Section \ref{sec:ER}, we show some numerical experiments which verify our theoretical analysis and demonstrate that, in comparison
with some often used conventional least-squares solvers, faster convergence  has been obtained by our method.
Finally, we end this work with some conclusions in Section \ref{sec:conclusions}.

\section{Preparation}\label{sec:related-algor}
Using Petrov-Galerkin (PG) conditions, Saad gave a prototype projection iteration format as follows \cite[Section 5.1.2]{03Saad}.
\begin{align}\label{eq:PG-iterate}
\bx^{(k+1)}  = \bx^{(k)} + V(W^T A V)^\dag W^T (\bb-A\bx^{(k)})
\end{align}
for $k=0,1,2,\cdots$, where $V$ and $W$ are two parameter matrices. By varying $V$ and $W$, one can obtain many popular iterations as special cases, including the following multiple row iterate scheme.

Let $W=\bet_k$ and $V=A^TW$, where $\bet_k\in\Rc^{m}$ is any non-zero vector. The iteration step in \eqref{eq:PG-iterate} reduces to
\begin{align}\label{eq:GK}
  \bx^{(k+1)}  = \bx^{(k)} + \frac{\bet_k^T(\bb - A\bx^{(k)})}{\BT{A^T \bet_k}} A^T \bet_k.
\end{align}
In particular, if $\bet_k=[\eta_{k,i}]_{i=1}^m\in \Rc^{m}$ is a Gaussian vector with $\eta_{k,i}\sim \Nc(0,1)$, Gower and Richt\'{a}rik \cite{15GR} called  \eqref{eq:GK} the Gaussian Kaczmarz (GK) method. From the fact that $\bet_k^T A =\sum_{i=1}^m \eta_{k,i} A_{i,:}^T$ is a linear combination of rows of $A$, we know that GK needs operate all of them at each iteration. Very recently, in Chen and Huang's work \cite{22CH1}, they proposed a fast deterministic block Kaczmarz method, in which a set $U_k$ is first computed according to the greedy index selection strategy \cite{18BW1} and then the vector $\bet_k$ is constructed by
$
  \bet_k = \sum_{i\in U_k} (b_i - A_{i,:} \bx^{(k)}) \bmu_{i}.
$
It follows that only partial rows indexed by $U_k$ are used. Its relaxed version was given by \cite{22WCZ}. However, this way may suffer from the time-consuming calculation. This is because, to compute $U_k$, one has to scan the residual vector from scratch during each iteration. It is unfavorable when the size of $A$ is large.  For more variants and studies of row iterative methods, we refer to \cite{21BW, 22BLZL, 15Dumitrescu, 22HNRS,  22SLTW, 212WX,17XZ} and the references therein.

Randomization has several benefits, e.g., the resulting algorithm is easy to analyze, simple to implement, and often effective in practice \cite{18BW1,15GR,14NT}. Here, we combine both ideas of randomization technique and  multiple row iterative schemes and propose a randomized multiple row (RMR) method and its variant respectively described in Algorithms \ref{alg:RMR} and \ref{alg:RMR2}.

\begin{algorithm}[htb]
\caption{ The \textbf{RMR} method for $A \bx= \bb$.}
\label{alg:RMR}
\begin{algorithmic}[1]
\Require
partition $\left\{\mathcal{I}_i\right\}_{i=1}^{s}$,
initial vector $\bx^{(0)}\in \Rc^{n}$, and maximum iteration number $\ell$.
\Ensure
$\bx^{(\ell)}$.
\State {\bf for} $k=0,1,2,\cdots,\ell$ {\bf do}
\State \quad pick $i\in[s]$ with probability $p_i = \frac{\BF{A_{\mathcal{I}_i,:}}}{\BF{A}}$;
\State \quad compute
\begin{align}\label{eq:RMR}
 \bx^{(k+1)}  = \bx^{(k)} + \frac{\bet_k^T(\bb -  A\bx^{(k)})}{\BT{A^T \bet_k}} A^T \bet_k,
\end{align}
       \quad where $\bet_k = \sum_{i\in \mathcal{I}_i} (b_i - A_{i,:} \bx^{(k)}) \bmu_{i}$;
\State {\bf endfor}
\end{algorithmic}
\end{algorithm}

\begin{algorithm}[htb]
\caption{ The \textbf{RMR} method for $A^T\by={\bm 0}$.}
\label{alg:RMR2}
\begin{algorithmic}[1]
\Require
partition $\left\{\mathcal{J}_j\right\}_{j=1}^{t}$,
initial vector $\by^{(0)}\in\Rc^{m}$, and maximum iteration number $\ell$.
\Ensure
$\by^{(\ell)}$.
\State {\bf for} $k=0,1,2,\cdots,\ell$ {\bf do}
\State \quad pick $j\in[t]$ with probability $q_j = \frac{\BF{A_{:,\mathcal{J}_j}}}{\BF{A}}$;
\State \quad compute
\begin{align}\label{eq:RMR+ATy=0}
 \by^{(k+1)}  = \by^{(k)} - \frac{\bze_k^T A^T \by^{(k)}}{\BT{A \bze_k}} A \bze_k,
\end{align}
       \quad where $\bze_k = \sum_{j\in \mathcal{J}_j} ( - A_{:,j}^T \by^{(k)}) \bnu_{j}$;
\State {\bf endfor}
\end{algorithmic}
\end{algorithm}

\vskip 0.5ex
Let $\bE\left[ \cdot \right]$ denote the expectation taken over the random choice of the algorithm. The convergence of RMR is characterized
by two rates $\rho_1$ and $\rho_2 \in [0,1)$ that depend on the smallest and largest singular values of a matrix respectively  defined by
\begin{align}\label{eq:two_para_rho}
  \rho_1    = 1 - \frac{1}{\beta_{\max}^{I}} \frac{\sigma_{\min}^2(A)}{\BF{A}}
  ~~{\rm and}~~
  \rho_2 = 1 - \frac{1}{\beta_{\max}^{J}} \frac{\sigma_{\min}^2(A)}{\BF{A}},
\end{align}
where
$\beta_{\max}^{I} := \max_{i\in [s]} \left\{ \sigma_{\max}^2(A_{\mathcal{I}_i,:})/\BF{A_{\mathcal{I}_i,:}}\right\}$
 and
$\beta_{\max}^{J} := \max_{j\in [t]} \left\{ \sigma_{\max}^2(A_{:,\mathcal{J}_j})/\BF{A_{:,\mathcal{J}_j}}\right\}$.

\vskip 0.5ex
\begin{theorem}\label{LemmaRMR}
  For any given consistent linear system $A\bx=\bb$, let $\bx^{(k)}$ be the vector generated by Algorithm \ref{alg:RMR} with $\bx^{(0)}$ being in the space $\Rc(A^T)$. It holds that
  \begin{align*}
    \bE\left[\BT{\bx^{(k)}- \bx^{\ast}} \right] \leq \rho_1^{k}~ \BT{\bx^{(0)}- \bx^{\ast}},
  \end{align*}
where the constant $\rho_1$ is defined by formula \eqref{eq:two_para_rho}.
\end{theorem}

\vskip 0.5ex
Instead of projecting directly onto the range of $A$ to find $\bx^{\ast}$, Algorithm \ref{alg:RMR} considers iterative projection on the range of smaller sketched matrix $\bet_k^TA$, where the random sketching vector $\bet_k$ is sampled from the partition of $[m]$. This projection-based method satisfies the following PG conditions. Specifically,
 \begin{align*}
 \bx^{(k+1)} \in \bx^{(k)} + {\rm span}\{A^T \bet_k\}
 ~~{\rm and}~~
 \bb - A\bx^{(k+1)} ~\bot~ {\rm span} \{\bet_k\}.
 \end{align*}
This problem can be solved directly and is equivalent to the iterative step in \eqref{eq:RMR}. Theorem \ref{LemmaRMR} shows that the approximation $\bx^{(k+1)}$  converges linearly in the mean square to the solution $\bx^{\ast}$ when the linear system is consistent. Details of proving this theorem are in the Appendix.

An interesting generalization is to use Algorithm \ref{alg:RMR} to solve a homogeneous linear system of equations $A^T\by  = {\bm 0}$, as shown in Algorithm \ref{alg:RMR2}. This system is underdetermined and consistent with infinitely many solutions. It has been shown that the orthogonal projection of $\bb$ onto $\Nc(A^T)$, equals $\bb_{N}$,  is one of the least-squares solutions with minimum Euclidean norm; see, e.g., \cite{211WX}. In the following, we give a convergence result of the algorithm.

\vskip 0.5ex

\begin{theorem}\label{ERMR:thm+RMR+Ay=0}
  For any given homogeneous linear system $A^T\by={\bm 0}$, let $\by^{(k)}$ be the vector generated by Algorithm \ref{alg:RMR2} with $\by^{(0)}$ being in the affine space $\bb+\Rc(A)$. It holds that
  \begin{align*}
    \bE\left[\BT{\by^{(k)}-\bb_{N} } \right] \leq \rho_2^{k} ~\BT{\by^{(0)}-\bb_{N}},
  \end{align*}
  where the constant $\rho_2$ is defined by formula \eqref{eq:two_para_rho}.
\end{theorem}

\vskip 0.5ex
Theorem \ref{ERMR:thm+RMR+Ay=0} bounds the convergence rate of Algorithm \ref{alg:RMR2} in the mean square. Initially, it starts with $\by^{(0)}\in \bb+\Rc(A)$. At the $k$th iteration, the algorithm randomly selects the index $j$ from $[t]$ and updates $\by^{(k)}$ by projecting it onto the orthogonal complement of $\Rc(A)$. The claim is that randomly selecting the index with a probability implies that the algorithm converges to $\bb_{N}$ in the mean square. After $k$ iterations,  the algorithm outputs $\by^{(k)}$ and by orthogonality $\bb-\by^{(k)}$ serves as an approximation of $\bb - \bb_{N} = \bb_{R}$ when $k \rightarrow \infty$. We defer the proof of Theorem \ref{ERMR:thm+RMR+Ay=0} to the Appendix.

\section{The extended randomized multiple row method}\label{sec:ERMR}
The analyses and discussions in the previous section mainly focus on consistent systems. In practice, we usually have access only to a contaminated  right-hand side with $\bb_N\neq 0$. This type of problem arises in many real-world applications; see, e.g., \cite{14AH, 14DL, 18HJ, 18LMD}.

When system \eqref{eq:Ax=b} is not consistent, from the convergence analysis of RMR, we mention that
\begin{align*}
\bx^{(k+1)}- \bx^{\ast}
 = \left(I_n - \frac{A^T \bet_k \bet_k^T A }{\BT{A^T \bet_k}}\right) \left(\bx^{(k)}- \bx^{\ast}\right) +
       \frac{\bet_k^T \bb_N }{\BT{A^T \bet_k}} A^T \bet_k
\end{align*}
with $A\bx^{\ast}=\bb_R$. By orthogonality, it follows that
\begin{align*}
\bE\left[  \BT{\bx^{(k+1)}- \bx^{\ast}} \right]
 = \bE\left[  \BT{\bx^{(k)}- \bx^{\ast}}   \right] - \bE\left[  \frac{|\bet_k^T A \left(\bx^{(k)}- \bx^{\ast}\right)|^2}{\BT{A^T \bet_k}} \right] +
    \bE\left[  \frac{|\bet_k^T \bb_N |^2}{\BT{A^T \bet_k}} \right].
\end{align*}
This fact tells us that RMR may not converge to the least-squares solution due to the  existence of $\bb_N$. To resolve this problem, inspired by the extended row method in \cite{20DSS, 14NT,15NZZ, 211WX}, we provide the following extended variant of RMR (ERMR)  and analyze its convergence and computational complexity.

\subsection{The proposed method}
At the $k$th iterate, ERMR first carries out Algorithm \ref{alg:RMR2} on the homogeneous linear system $A^T \by  = {\bm 0}$ and output $\by^{(k+1)}$, and then applies Algorithm \ref{alg:RMR} to the asymptotical consistent linear system $A \bx   = \bb- \by^{(k+1)}\approx \bb_R$ and output $\bx^{(k+1)}$. The ERMR method is formally described in Algorithm \ref{alg:ERMR}.

\begin{algorithm}[htb]
\caption{ The \textbf{ERMR} method for $A \bx = \bb$.}
\label{alg:ERMR}
\begin{algorithmic}[1]
\Require
partitions $\left\{\mathcal{I}_i\right\}_{i=1}^{s}$ and $\left\{\mathcal{J}_j\right\}_{j=1}^{t}$,
initial vectors $\bx^{(0)}\in \Rc^{n}$ and $\by^{(0)}\in\Rc^{m}$, and maximum iteration number $\ell$.
\Ensure
$\bx^{(\ell)}$.
\State {\bf for} $k=0,1,2,\cdots,\ell$ {\bf do}
\State \quad pick $j\in[t]$ with probability $q_j = \frac{\BF{A_{:,\mathcal{J}_j}}}{\BF{A}}$;
\State \quad compute
\begin{align}\label{eq:ERMR1}
 \by^{(k+1)}  = \by^{(k)} - \frac{\bze_k^T A^T \by^{(k)}}{\BT{A \bze_k}} A \bze_k,
\end{align}
       \quad where $\bze_k = \sum_{j\in \mathcal{J}_j} ( - A_{:,j}^T \by^{(k)}) \bnu_{j}$;
\State \quad pick $i\in[s]$ with probability $p_i = \frac{\BF{A_{\mathcal{I}_i,:}}}{\BF{A}}$;
\State \quad compute
\begin{align}\label{eq:ERMR2}
 \bx^{(k+1)}  = \bx^{(k)} + \frac{\bet_k^T(\bb - \by^{(k+1)} - A\bx^{(k)})}{\BT{A^T \bet_k}} A^T \bet_k,
\end{align}
       \quad where $\bet_k = \sum_{i\in \mathcal{I}_i} (b_i -  y^{(k+1)}_i - A_{i,:} \bx^{(k)}) \bmu_{i}$;
\State {\bf endfor}
\end{algorithmic}
\end{algorithm}

In the following, we give two remarks to illustrate the relationships between ERMR and several existing extended randomized row iterative methods, including REABK \cite{20DSS} and GEK \cite{211WX}.

\vskip 0.5ex
\begin{remark}\label{remark:GEKvsERMR}
 The updates in our method are reminiscent of GEK \cite{211WX}, defined by \eqref{gek}. Indeed they are all the pseudoinverse-free REK methods. At the $k$th iterate, the quantity $ \bet^T A $ (resp., $\bze^T A^T$) in updating $\bx^{(k)}$ (resp., $\by^{(k)}$) is in the form of a linear combination of all rows of $A$ (resp., $A^T$). It enables GEK that accesses all rows of $A$ and $A^T$ each iteration. Nonetheless, in general, computing this way is very demanding in terms of computational efficiency. To compute the next update at each ERMR iteration step, our method only partially utilizes the rows concerning the randomized index sets $\mathcal{I}_i$ and $\mathcal{J}_j$.  This local format can save computational resources significantly. For more details, we refer to the computational complexity analysis of ERMR in Section \ref{ERMR+subsec:CC}.
\end{remark}

\vskip 0.5ex
\begin{remark}
In fact, each half REABK iterate, defined by \eqref{reabk}, is the projection of the previous iterate onto each row indexed by a set and then averaged. By direct computation, the updates \eqref{eq:ERMR1} and \eqref{eq:ERMR2} respectively become
\begin{align*}
 \by^{(k+1)}
 & = \by^{(k)} -
 \frac
 { \sum_{j\in \mathcal{J}_j} ( - A_{:,j}^T \by^{(k)}) \bnu_{j}^T A^T \by^{(k)} }
 { \BT{A \bze_k} }
 \sum_{j\in \mathcal{J}_j} ( - A_{:,j}^T \by^{(k)}) A\bnu_{j}\\
 & = \by^{(k)} - \frac
 { \BT{\bze_k}  }{ \BT{A \bze_k} } A_{:,\mathcal{J}_j}A_{:,\mathcal{J}_j}^T \by^{(k)}\\
 & = \by^{(k)} - \widehat{\alpha}_k\frac{A_{:,\mathcal{J}_j}A_{:,\mathcal{J}_j}^T \by^{(k)}}{\BF{A_{:,\mathcal{J}_j}}}
\end{align*}
and
\begin{align*}
 \bx^{(k+1)}
 & = \bx^{(k)} +
 \frac{\BT{\bet_k} }{\BT{A^T \bet_k}} A_{\mathcal{I}_i,:}^T \left(\bb_{\mathcal{I}_i}- \by^{(k)}_{\mathcal{I}_i} - A_{\mathcal{I}_i,:} \bx^{(k)} \right)\\
 & = \bx^{(k)} + \widetilde{\alpha}_k\frac{A_{\mathcal{I}_i,:}^T \left(\bb_{\mathcal{I}_i}- \by^{(k)}_{\mathcal{I}_i} - A_{\mathcal{I}_i,:} \bx^{(k)} \right)}{\BF{A_{\mathcal{I}_i,:}}},
\end{align*}
where the  adaptive step sizes $\widehat{\alpha}_k = \BT{\bze_k} \BF{A_{:,\mathcal{J}_j}}/\BT{A \bze_k}$ and $\widetilde{\alpha}_k = \BT{\bet_k} \BF{A_{\mathcal{I}_i,:}}/\BT{A^T\bet_k}$.
Then, ERMR can be viewed as a asynchronous REABK in which the latter uses a constant step size, see \cite[Algorithm 1]{20DSS}.
\end{remark}

\subsection{Convergence analysis}
 In this section, we will prove that ERMR lineally converges in the mean square to $\bx^{\ast}$ when the system \eqref{eq:Ax=b} is not consistent.

\vskip 0.5ex

\begin{theorem}\label{ERMR:TheoremERMR}
  For any given inconsistent linear system $A\bx=\bb$, let $\bx^{(k)}$ be the vector generated by Algorithm \ref{alg:ERMR} with
  $\bx^{(0)}$ being in the space $\Rc(A^T)$ and $\by^{(0)}$ being in the affine space $\bb+\Rc(A)$. It holds that
  \begin{align}\label{ERMR:TheoremResult1}
    \bE\left[\BT{\bx^{(k)} - \bx^{\ast}} \right] \leq   \omega\rho^{k},
  \end{align}
 where $\rho = \max\{ \rho_1 ,   \rho_2\}$ $(\rho_1  \neq \rho_2)$ with $\rho_1$ and $\rho_2$ being defined by formula \eqref{eq:two_para_rho} and $\omega$ is defined by
 \begin{align*}
   \omega:= \BT{\bx^{(0)} - \bx^{\ast}} + \frac{s}{|\rho_1 - \rho_2|} \frac{\BT{\by^{(0)}-\bb_{N}}}{\BF{A}\beta_{\min}^{I}},
 \end{align*}
 with $\beta_{\min}^{I} = \min_{i\in [s]} \left\{ \sigma_{\min}^2(A_{\mathcal{I}_i,:})/\BF{A_{\mathcal{I}_i,:}}\right\}$.
\end{theorem}

\vskip 0.5ex
\begin{proof}
For $k=0,1,2,\cdots$, define an auxiliary vector
 \begin{align*}
 \widetilde{\bx}^{(k+1)}  = \bx^{(k)} + \frac{\bet_k^T\left(\bb_{R} - A\bx^{(k)}\right)}{\BT{A^T \bet_k}} A^T \bet_k.
\end{align*}
The error $\bx^{(k+1)} - \bx^{\ast}$ consists of $\bx^{(k+1)} - \widetilde{\bx}^{(k+1)}$ and $\widetilde{\bx}^{(k+1)}- \bx^{\ast}$. More specifically,
\begin{align*}
\bx^{(k+1)} - \widetilde{\bx}^{(k+1)}
= \frac{\bet_k^T\left(\bb_{N} - \by^{(k+1)}\right)}{\BT{A^T \bet_k}} A^T \bet_k
\end{align*}
and
\begin{align*}
\widetilde{\bx}^{(k+1)}- \bx^{\ast} = \left(I_n - \frac{A^T \bet_k \bet_k^TA}{\BT{A^T \bet_k}} \right)\left( \bx^{(k)} - \bx^{\ast} \right).
\end{align*}
By the orthogonality, namely,
\begin{align*}
\left(A^T \bet_k\right)^T \left(I_n - \frac{A^T \bet_k \bet_k^TA}{\BT{A^T \bet_k}} \right)\left( \bx^{(k)} - \bx^{\ast} \right)=0,
\end{align*}
the expected squared error admits
\begin{align}\label{thm1:eq0}
  \bE\left[\BT{\bx^{(k+1)} - \bx^{\ast}}\right]
  = \underbrace{\bE\left[\BT{\bx^{(k+1)} - \widetilde{\bx}^{(k+1)}}\right]}_{{\rm \textbf{Term 1}}} +
    \underbrace{\bE\left[\BT{\widetilde{\bx}^{(k+1)}- \bx^{\ast}}\right]}_{{\rm \textbf{Term 2}}}
\end{align}
with
\begin{align}\label{thm1:eq1}
\BT{\bx^{(k+1)} - \widetilde{\bx}^{(k+1)}}
= \frac{|\bet_k^T\left(\bb_{N} - \by^{(k+1)}\right)|^2}{\BT{A^T \bet_k}}
\end{align}
and
\begin{align}\label{thm1:eq2}
\BT{\widetilde{\bx}^{(k+1)}- \bx^{\ast}}
= \BT{\bx^{(k)} - \bx^{\ast}} - \frac{|\bet_k^T A \left( \bx^{(k)} - \bx^{\ast} \right)|^2 }{\BT{A^T \bet_k}}.
\end{align}
We proceed to analyze them individually.

{\textbf{Term 1}:} As already known, we use the same notation in \cite{20DSS, 13ZF}. The conditional expectation conditioned on the first $k$ iterations is defined by
\begin{align*}
\bE_{k}[\cdot] =\bE [\cdot|j_0,i_0,j_1,i_1,\cdots,j_{k},i_{k}],
\end{align*}
 where  $i_{\ell}$ and $j_\ell$ mean that the index sets $\left\{\mathcal{I}_i\right\}_{i=0}^{s}$ and $\left\{\mathcal{J}_j\right\}_{j=0}^{t}$ are chosen as
\begin{align*}
 \bE_{k}^j\left[\cdot\right]=\bE \left[\cdot|j_0,i_0 ,\cdots,j_{k},i_{k},i_{k+1}\right]
 \quad {\rm and} \quad
 \bE_{k}^i\left[\cdot\right]=\bE \left[\cdot|j_0,i_0,\cdots,j_{k},i_{k},j_{k+1}\right],
\end{align*}
respectively.  Since the random variables $i$ and $j$ are independent of each other, by the law of total expectation, we have
\begin{equation*}
  \bE_{k}\left[\cdot\right]=\bE_{k}^j\left[\bE_{k}^i\left[\cdot\right]\right]
  ~{\rm and}~
  \bE_{k}\left[\cdot\right]=\bE_{k}^i\left[\bE_{k}^j\left[\cdot\right]\right].
\end{equation*}
An elementary computation shows that
\begin{align*}
\frac{\BF{A_{\mathcal{I}_i,:}}}{\BT{A^T \bet_k}}
& = \frac{\BF{A_{\mathcal{I}_i,:}}}{\BT{A_{\mathcal{I}_i,:}^T\left( \bb_{\mathcal{I}_i} - \by_{\mathcal{I}_i}^{(k+1)}- A_{\mathcal{I}_i,:}\bx^{(k)} \right)}}\\
& \leq \frac{\BF{A_{\mathcal{I}_i,:}}}{\sigma_{\min}^2\left( A_{\mathcal{I}_i,:}\right) \BT{ \bb_{\mathcal{I}_i} - \by_{\mathcal{I}_i}^{(k+1)}- A_{\mathcal{I}_i,:}\bx^{(k)} }}\\
& \leq \frac{1}{\beta_{\min}^I} \frac{1}{\BT{ \bb_{\mathcal{I}_i} - \by_{\mathcal{I}_i}^{(k+1)}- A_{\mathcal{I}_i,:}\bx^{(k)} }}
= \frac{1}{\beta_{\min}^I}\frac{1}{\BT{ \bet_k }},
\end{align*}
where the second line is from the fact that $\bb_{\mathcal{I}_i} - \by_{\mathcal{I}_i}^{(k+1)}- A_{\mathcal{I}_i,:}\bx^{(k)} \in \Rc(A_{\mathcal{I}_i,:})$, it follows that
\begin{align*}
\frac{|\bet_k^T\left(\bb_{N} - \by^{(k+1)}\right)|^2}{\BT{A^T \bet_k}}
& \leq \frac{\BT{\bet_k} \BT{\by^{(k+1)} - \bb_{N}}}{\BT{A^T \bet_k}}
\leq \frac{1}{\beta_{\min}^I} \frac{\BT{\by^{(k+1)} - \bb_{N}}}{\BF{A_{\mathcal{I}_i,:}}}.
\end{align*}
Then, we have
\begin{align*}
\bE_k\left[\BT{\bx^{(k+1)} - \widetilde{\bx}^{(k+1)}}\right]
&\leq \frac{1}{\beta_{\min}^I} \bE_{k}^j\left[\bE_{k}^i\left[ \frac{\BT{\by^{(k+1)} - \bb_{N}}}{\BF{A_{\mathcal{I}_i,:}}}  \right]\right] \\
& = \frac{1}{\beta_{\min}^I} \bE_{k}^j\left[ \sum_{i=1}^s \frac{\BF{A_{\mathcal{I}_i,:}}}{\BF{A}} \frac{\BT{\by^{(k+1)} - \bb_{N}}}{\BF{A_{\mathcal{I}_i,:}}}  \right] \\
& = \frac{s}{\BF{A}\beta_{\min}^I} \bE_{k}^j \left[ \BT{\by^{(k+1)} - \bb_{N}} \right].
\end{align*}
By taking expectation again, Theorem \ref{ERMR:thm+RMR+Ay=0} gives
\begin{align}\label{thm1:eq4}
\bE\left[\BT{\bx^{(k+1)} - \widetilde{\bx}^{(k+1)}}\right]
\leq \frac{s}{\BF{A}\beta_{\min}^I} \bE \left[ \BT{\by^{(k+1)} - \bb_{N}} \right]
\leq \frac{s\rho_2^{k+1}}{\BF{A}\beta_{\min}^I} \BT{\by^{(0)} - \bb_{N}}.
\end{align}

{\textbf{Term 2}:} Based on the following two facts that
\begin{align*}
  \left|\bet_k^T A \left( \bx^{(k)} - \bx^{\ast} \right)\right|^2
  & = \left|\sum_{i\in \mathcal{I}_i}\left(\bb_i -\by_i^{(k+1)} - A_{i,:}\bx^{(k)}\right) \bmu_i^T A \left( \bx^{(k)} - \bx^{\ast} \right)\right|^2\\
  & = \BT{\left( \bb_{\mathcal{I}_i} - \by_{\mathcal{I}_i}^{(k+1)}- A_{\mathcal{I}_i,:}\bx^{(k)} \right)^T A_{\mathcal{I}_i,:}\left( \bx^{(k)} - \bx^{\ast} \right)}\\
  & \leq \BT{\bb_{\mathcal{I}_i} - \by_{\mathcal{I}_i}^{(k+1)}- A_{\mathcal{I}_i,:}\bx^{(k)}} \BT{ A_{\mathcal{I}_i,:}\left( \bx^{(k)} - \bx^{\ast} \right)}
\end{align*}
and
\begin{align*}
\frac{\BF{A_{\mathcal{I}_i,:}}}{\BT{A^T \bet_k}}
& \geq \frac{\BF{A_{\mathcal{I}_i,:}}}{\sigma_{\max}^2\left( A_{\mathcal{I}_i,:}\right) \BT{ \bb_{\mathcal{I}_i} - \by_{\mathcal{I}_i}^{(k+1)}- A_{\mathcal{I}_i,:}\bx^{(k)} }}\\
& \geq \frac{1}{\beta_{\max}^I} \frac{1}{\BT{ \bb_{\mathcal{I}_i} - \by_{\mathcal{I}_i}^{(k+1)}- A_{\mathcal{I}_i,:}\bx^{(k)} }},
\end{align*}
the conditioned expectation for the second term on the right of formula \eqref{thm1:eq2} is bounded by
\begin{align*}
\bE_k^i\left[-\frac{\big|\bet_k^T A \left( \bx^{(k)} - \bx^{\ast} \right)\big|^2 }{\BT{A^T \bet_k}}\right]
& = \bE_k^i
\left[
-\frac{\BF{A_{\mathcal{I}_i,:}}}{\BT{A^T \bet_k}} \cdot
\big|\bet_k^T A \left( \bx^{(k)} - \bx^{\ast} \right)\big|^2 \cdot
\frac{1}{\BF{A_{\mathcal{I}_i,:}}}
\right]\\
& \leq -\frac{1}{\beta_{\max}^{I}} \bE_k^i\left[\frac{\BT{A_{\mathcal{I}_i,:}\left(\bx^{(k)} - \bx^{\ast}\right)}}{\BF{A_{\mathcal{I}_i,:}}} \right]\\
& = -\frac{1}{\beta_{\max}^{I}}
\sum_{i=1}^s \frac{\BF{A_{\mathcal{I}_i,:}}}{\BF{A}} \frac{\BT{A_{\mathcal{I}_i,:}\left(\bx^{(k)} - \bx^{\ast}\right)}}{\BF{A_{\mathcal{I}_i,:}}}\\
& = -\frac{\BT{A\left(\bx^{(k)} - \bx^{\ast}\right)}}{\beta_{\max}^{I} \BF{A}}\\
& \leq -\frac{1}{\beta_{\max}^{I}} \frac{\sigma_{\min}^2(A)}{\BF{A}}\BT{\bx^{(k)} - \bx^{\ast}}.
  \end{align*}
By taking expectation for both sides of the above inequality, it holds that
  \begin{equation}\label{thm1:eq3}
    \bE\left[\BT{\widetilde{\bx}^{(k+1)}- \bx^{\ast}}\right]
    \leq \rho_1 ~ \bE\left[\BT{\bx^{(k)} - \bx^{\ast}} \right].
  \end{equation}

Combining formulas \eqref{thm1:eq0}, \eqref{thm1:eq4}, and \eqref{thm1:eq3}, it indicates that
\begin{align*}
  \bE\left[\BT{\bx^{(k+1)} - \bx^{\ast}}\right]
  \leq \frac{s\rho_2^{k+1}}{\BF{A}\beta_{\min}^I} \BT{\by^{(0)} - \bb_{N}} + \rho_1 ~ \bE\left[\BT{\bx^{(k)} - \bx^{\ast}} \right].
\end{align*}
By unrolling the recurrence, we obtain the estimate
\begin{align*}
\bE\left[\BT{\bx^{(k)} - \bx^{\ast}} \right] \leq
    \rho_1^{k} ~ \BT{\bx^{(0)} - \bx^{\ast}} +
    \frac{s}{\BF{A}\beta_{\min}^{I}} \sum_{\ell=0}^{k-1}\rho_2^{\ell}\rho_1^{k-\ell}~ \BT{\by^{(0)}-\bb_{N}}.
\end{align*}
 Note that if $\rho_1  < \rho_2$, we have
\begin{equation*}
 \sum_{\ell=0}^{k-1}\rho_2^{k-\ell} \rho_1^\ell
 = \rho_2^{k} \sum_{j=0}^{k-1} \left( \frac{\rho_1}{\rho_2}\right)^\ell
 \leq \frac{\rho_2^{k}}{\rho_2 -\rho_1}
 \leq  \frac{\rho^{k}}{|\rho_1  - \rho_2|}.
\end{equation*}
If $\rho_2 <\rho_1$,  the result also holds. It follows that
\begin{align*}
\bE\left[\BT{\bx^{(k)} - \bx^{\ast}} \right]
      \leq \rho^{k} ~ \BT{\bx^{(0)} - \bx^{\ast}} +
           \rho^{k} ~ \frac{s}{|\rho_1  - \rho_2|\BF{A}\beta_{\min}^{I}} \BT{\by^{(0)}-\bb_{N}} \leq   \omega\rho^{k}.
\end{align*}
This concludes the proof of the theorem.
\hfill
\end{proof}

\vskip .5ex
\begin{remark}
By iteratively removing the component in $\bb_{N}$, the current vector $\bx^{(k)}$ converges to $\bx^{\ast}$. Step $3$ in Algorithm \ref{alg:ERMR} plays an important role in ERMR. But, when the linear system is consistent, we do not need the iteration of $\by^{(k)}$ anymore. In this case, we set $\by^{(k)}=0$ for $k=0,1,2,\cdots$ in ERMR and compute $\bx^{(k)}$ according to the RMR method in Algorithm \ref{alg:RMR}.
\end{remark}

\subsection{Computational complexity}\label{ERMR+subsec:CC}
In this section we will discuss the computational complexity of Algorithm \ref{alg:ERMR} as follows.

\vskip 0.5ex

\begin{theorem}\label{ERMR:lemma1}
For any given inconsistent linear system $A\bx=\bb$, let $\bx^{(k)}$ be the vector generated by Algorithm \ref{alg:ERMR} with $\bx^{(0)}$ being in the space $\Rc(A^T)$ and $\by^{(0)}$ being in the affine space $\bb+\Rc(A)$.  Given an accuracy tolerance $\epsilon$  and a confidence level $\beta$, when the iteration number satisfies $k \geq (1-\rho)^{-1} \ln( \omega /(\varepsilon(1-\beta))$,
it holds that
\begin{equation}\label{ERMR-Lemma3:CC}
\bP\left[\BT{\bx^{(k)} - \bx^{\ast}}  < \varepsilon \right] \geq \beta,
\end{equation}
where the constants $\omega$ and $\rho$ are defined by Theorem \ref{ERMR:TheoremERMR}.
\end{theorem}

\vskip 0.5ex

\begin{proof}
When $k \geq (1-\rho)^{-1}\ln( \omega /(\varepsilon(1-\beta))$, we can check that for $0<\rho<1$,
\begin{equation*}
 \rho^k
 \leq e^{-(1-\rho)k}
 \leq e^{-\ln\frac{\omega}{\varepsilon(1-\beta)}}
 =\frac{\varepsilon(1-\beta)}{\omega}.
\end{equation*}
Using Markov's inequality,  which says that for a positive random variable $x$ and an accuracy parameter $\varepsilon\in(0,1)$, the probability that $x$ is greater than or equal to $\varepsilon$  is less than or equal to the expected value of $x$ divided by $\varepsilon$, the claim is concluded as follows.
\begin{align*}
     \bP\left[\BT{\bx^{(k)} - \bx^{\ast}}  \geq \varepsilon \right]
       \leq \frac{1}{\varepsilon} \bE\left[\BT{\bx^{(k)} - \bx^{\ast}}\right]
       \leq \frac{\omega}{\varepsilon} \rho^{k}
       \leq 1-\beta,
  \end{align*}
where the second inequality is from  formula \eqref{ERMR:TheoremResult1}.
\hfill
\end{proof}

\vskip 0.5ex

When the probability of the squared error $\BT{\bx^{(k)} - \bx^{\ast}}$ below $\varepsilon$  is at least  $\beta$, Theorem \ref{ERMR:lemma1} bounds the number of iterations required by the algorithm to terminate. After that, we analyze the flopping operations (flops) per iteration. Let $\widetilde{A}=AA^T$ and $\widehat{A}=A^TA$. The pre-processing of  computing them is dominant in flops, but these can be very efficiently performed by the high-level basic linear algebra subroutine, e.g., BLAS3.  Thus, we always assume the availability of $\widetilde{A}$ and $\widehat{A}$.

At step $3$ in Algorithm \ref{alg:ERMR}, define an additional vector $\widehat{\br}^{(k)}=-A^T\by^{(k)}$ for $k=0,1,2,\cdots$. We have the following recursive formula.
\begin{align*}
 \widehat{\br}^{(k+1)}
 & = \widehat{\br}^{(k)} -
 \frac{\sum_{j\in \mathcal{J}_j} \widehat{r}^{(k)}_{j}\bnu_{j}^T \widehat{\br}^{(k)}}
      {\sum_{j\in \mathcal{J}_j} \sum_{h\in \mathcal{J}_j} \widehat{r}^{(k)}_{j}  \widehat{r}^{(k)}_{h} \bnu_{j} A^T A \bnu_{h}^T}
       \sum_{j\in \mathcal{J}_j}  \widehat{r}^{(k)}_{j} A^T A \bnu_{j}\\
 & = \widehat{\br}^{(k)} -
 \frac{\sum_{j\in \mathcal{J}_j}\left( \widehat{r}^{(k)}_{j}\right)^2}
      {\sum_{j\in \mathcal{J}_j} \sum_{h \in \mathcal{J}_j}  \widehat{r}^{(k)}_{j} \widehat{r}^{(k)}_{h} \widehat{A}_{j, h} }
       \sum_{j\in \mathcal{J}_j}  \widehat{r}^{(k)}_{j} \widehat{A}_{:,j}.
\end{align*}
As for $\by^{(k)}$, there is a recursion
\begin{align*}
 \by^{(k+1)}
  = \by^{(k)} -
 \frac{\sum_{j\in \mathcal{J}_j}\left( \widehat{r}^{(k)}_{j}\right)^2}
      {\sum_{j\in \mathcal{J}_j} \sum_{h \in \mathcal{J}_j}  \widehat{r}^{(k)}_{j} \widehat{r}^{(k)}_{h} \widehat{A}_{j, h} }
       \sum_{j\in \mathcal{J}_j}  \widehat{r}^{(k)}_{j} A_{:,j}.
\end{align*}

At step $5$ in Algorithm \ref{alg:ERMR}, let $\widetilde{\br}^{(k)}=\bb -\by^{(k+1)}-A\bx^{(k)}$ on $\bx^{(k)}$ be another supplementary vector for $k=0,1,2,\cdots$. It follows that
\begin{align*}
 \widetilde{\br}^{(k+1)}
 & = \bb -\by^{(k+1)}-A\bx^{(k+1)}\\
 & = \widetilde{\br}^{(k)} -
 \frac{\sum_{i\in \mathcal{I}_i} \widetilde{r}^{(k)}_{i}\bmu_{i}^T \widetilde{\br}^{(k)}}
      {\sum_{i\in \mathcal{I}_i} \sum_{g\in \mathcal{I}_i} \widetilde{r}^{(k)}_{i}  \widetilde{r}^{(k)}_{g} \bmu_{i} AA^T \bmu_{g}^T}
       \sum_{i\in \mathcal{I}_i}  \widetilde{r}^{(k)}_{i} AA^T \bmu_{i}\\
 & = \widetilde{\br}^{(k)} -
 \frac{\sum_{i\in \mathcal{I}_i}\left( \widetilde{r}^{(k)}_{i}\right)^2}
      {\sum_{i\in \mathcal{I}_i} \sum_{g\in \mathcal{I}_i}  \widetilde{r}^{(k)}_{i} \widetilde{r}^{(k)}_{g} \widetilde{A}_{i, g} }
       \sum_{i\in \mathcal{I}_i}  \widetilde{r}^{(k)}_{i} \widetilde{A}_{:,i}.
\end{align*}
Then, the  next approximation is computed by
\begin{align*}
 \bx^{(k+1)}
  = \bx^{(k)} +
 \frac{\sum_{i\in \mathcal{I}_i}\left( \widetilde{r}^{(k)}_{i}\right)^2}
      {\sum_{i\in \mathcal{I}_i} \sum_{g\in \mathcal{I}_i}  \widetilde{r}^{(k)}_{i} \widetilde{r}^{(k)}_{g} \widetilde{A}_{i, g} }
       \sum_{i\in \mathcal{I}_i}  \widetilde{r}^{(k)}_{i} A_{i,:}^T.
\end{align*}

Based on the recursive update formulas of $\by^{(k)}$ and $\bx^{(k)}$, Algorithm \ref{alg:ERMR} can be implemented efficiently because we do not have to directly compute $A \bze_k$ and $A^T \bet_k$. We emphasize that the algorithm updates only two vector $\by^{(k)}$ and $\bx^{(k)}$ (and uses two intermediate vectors $\widehat{\br}^{(k)}$ and $\widetilde{\br}^{(k)}$). We always apply blocking strategy to the matrix-vector multiplication, e.g.,
\begin{align*}
   \sum_{j\in \mathcal{J}_j}  \widehat{r}^{(k)}_{j} \widehat{A}_{:,j} = \widehat{A}_{:,\mathcal{J}_j}\widehat{\br}^{(k)}_{\mathcal{J}_j}
  ~~{\rm and}~~
   \sum_{j\in \mathcal{J}_j} \sum_{h \in \mathcal{J}_j}  \widehat{r}^{(k)}_{j} \widehat{r}^{(k)}_{h} \widehat{A}_{j, h}
 =\left( \widehat{\br}^{(k)}_{\mathcal{J}_j}\right)^T
  \left( \widehat{A}_{\mathcal{J}_j, \mathcal{J}_j} \widehat{\br}^{(k)}_{\mathcal{J}_j}\right).
\end{align*}
The computational processes are summarized in Table \ref{ERMR-tab:ComputXk}. Totally, it requires $\mathcal{O}(m)$
 flops.

\begin{table}[!ht]
\centering\renewcommand\arraystretch{1.25}
    \caption{The complexity of computing $\bx^{(k+1)}$ in ERMR.}
    \begin{tabular}{p{2cm}cp{6cm}cp{3cm}}
    \hline
    \multicolumn{3}{l}{{\sf Computing}~ $\by^{(k+1)}$}  \\
    \hline
   {\sf Step~1} && $g_1 = \left( \widehat{\br}^{(k)}_{\mathcal{J}_j}\right)^T \left( \widehat{\br}^{(k)}_{\mathcal{J}_j}\right)$  && $2|\mathcal{J}_j|-1$ \\
   {\sf Step~2} &&
   $g_2 =
   \left( \widehat{\br}^{(k)}_{\mathcal{J}_j}\right)^T
   \left( \widehat{A}_{\mathcal{J}_j, \mathcal{J}_j} \widehat{\br}^{(k)}_{\mathcal{J}_j}\right)
  $  &&$ 3(|\mathcal{J}_j|^2+|\mathcal{J}_j|)/2 $  \\
   {\sf Step~3} && $g_3 = g_1/g_2$                                 && $1$  \\
   {\sf Step~4} && ${\bm g}_1 = \widehat{A}_{:,\mathcal{J}_j}\widehat{\br}^{(k)}_{\mathcal{J}_j}$   && $n(2|\mathcal{J}_j|-1)$  \\
   {\sf Step~5} && $\widehat{\br}^{(k+1)} = \widehat{\br}^{(k)} - g_3 \cdot {\bm g}_1$  && $2n$  \\
   {\sf Step~6} && ${\bm g}_2 = A_{:,\mathcal{J}_j}\widehat{\br}^{(k)}_{\mathcal{J}_j}$  && $m(2|\mathcal{J}_j|-1)$  \\
   {\sf Step~7} && $\by^{(k+1)} = \by^{(k)} - g_3 \cdot {\bm g}_2$             && $2m$  \\
   \hline
   \multicolumn{3}{l}{{\sf Computing}~ $\bx^{(k+1)}$}  \\
   \hline
   {\sf Step~1} && $g_4 = \left( \widetilde{\br}^{(k)}_{\mathcal{I}_i}\right)^T \left( \widetilde{\br}^{(k)}_{\mathcal{I}_i}\right)$  && $2|\mathcal{I}_i|-1$ \\
   {\sf Step~2} && $g_5 = \left( \widetilde{\br}^{(k)}_{\mathcal{I}_i}\right)^T
   \left( \widetilde{A}_{\mathcal{I}_i, \mathcal{I}_i} \widetilde{\br}^{(k)}_{\mathcal{I}_i} \right)$  &&$ 3(|\mathcal{I}_i|^2+|\mathcal{I}_i|)/2 $  \\
   {\sf Step~3} && $g_6 = g_4/g_5$                                 && $1$  \\
   {\sf Step~4} && ${\bm g}_3 = \widetilde{A}_{:,\mathcal{I}_i} \widetilde{\br}^{(k)}_{\mathcal{I}_i}$   && $m(2|\mathcal{I}_i|-1)$  \\
   {\sf Step~5} && $\widetilde{\br}^{(k+1)} = \widetilde{\br}^{(k)} - g_6 \cdot {\bm g}_3$  && $2m$  \\
   {\sf Step~6} && ${\bm g}_4 = A_{\mathcal{I}_i,:}^T \widetilde{\br}^{(k)}_{\mathcal{I}_i}$  && $n(2|\mathcal{I}_i|-1)$  \\
   {\sf Step~7} && $\bx^{(k+1)} = \bx^{(k)} + g_6 \cdot {\bm g}_4$             && $2n$ \\
   \hline
   \end{tabular}
   \begin{tablenotes}
   \footnotesize
   \item[1.] {\it The symbol $|\mathcal{I}|$ denotes the cardinality of a set $\mathcal{I}$.}
   \end{tablenotes}
   \label{ERMR-tab:ComputXk}
\end{table}

\section{Experimental results}\label{sec:ER}
In this section we will give a few examples, including synthetic and real-world data, to demonstrate the convergence behaviors of RMR and ERMR.

In the spirit of randomization technique given by \cite{20DSS}, we suppose that the subsets $\left\{ \mathcal{I}_{i} \right\}_{i=1}^{s}$ and $\left\{ \mathcal{J}_{j} \right\}_{j=1}^{t}$ are respectively computed by
\begin{equation*}
\mathcal{I}_{i} =
 \left\{
 \begin{array}{ll}
 \left\{(i-1)\tau+1, (i-1)\tau+2,\cdots, i \tau \right\},& i\in [s-1], \vspace{1ex}\\
 \left\{(s-1)\tau+1, (s-1)\tau+2,\cdots, m      \right\},& i=s,
 \end{array}
 \right .
\end{equation*}
 and
 \begin{equation*}
\mathcal{J}_{j} =
 \left\{
 \begin{array}{ll}
 \left\{(j-1)\tau+1, (j-1)\tau+2,\cdots, j \tau \right\},& j\in [t-1], \vspace{1ex}\\
 \left\{(t-1)\tau+1, (t-1)\tau+2,\cdots, n      \right\},& j=t,
 \end{array}
 \right .
\end{equation*}
where $\tau$ is the size of $\mathcal{I}_{i}$ and $\mathcal{J}_{j}$. The discrete sampling is realized by applying MATLAB built-in function, e.g., {\sf  randsample}.

For comparisons, we use the implementations of REABK \cite[Algorithm 1]{20DSS} and GEK \cite{211WX}. As shown in \cite{20DSS}, we execute REABK using an empirical optimal step size
$\alpha       = 1.75/\beta_{\max}$  with
$\beta_{\max} = \max \left\{\beta_{\max}^{I}, \beta_{\max}^{J} \right\}$.
 Note that REABK is performed without explicitly forming $\alpha$. The algorithms are carried out on a Founder desktop PC with Intel(R) Core(TM) i5-7500 CPU 3.40 GHz.

\subsection{Synthetic data} The following coefficient matrix is generated from synthetic data.

\vskip 0.5ex
\begin{example}\label{ERMR:example1}
  As in Du et al. \cite{20DSS}, for given $m$, $n$, $r$, and $\kappa>1$, we construct a dense matrix $A$ by $A=UDV^T$, where $U\in \Rc^{m\times r}$, $D\in \Rc^{r\times r}$, and $V\in \Rc^{n\times r}$. Using MATLAB colon notation, these matrices are generated by
  $[U,\sim]={\sf qr}({\sf randn}(m,r),0)$,
  $[V,\sim]= {\sf qr}({\sf randn}(n,r),0)$, and
  $D={\sf diag}(1 + (\kappa - 1).* {\sf rand}(r,1))$.
\end{example}

\vskip 0.5ex
This example gives us many flexibilities to adjust the input parameters $m$, $n$, $r$, and $\kappa$ and yield various specific instantiations about the linear systems. In the following, we consider two types of {\it rank-deficient} cases by setting $m=30n$, $r=n/2$, and $\kappa=n/10$, which results in a mildly ill-posed test problem.

Since $A$ is rank-deficient, the right-hand side in \eqref{eq:Ax=b} is constructed by $\bb = A \bx^{\ast} + \delta \cdot \widehat{\bb}/\bT{\widehat{\bb}}$ with $\bx^{\ast}=A^{\dag} \widetilde{\bb}$, where $\widetilde{\bb}\in \Rc^m$ is generated from a standard normal distribution, $\widehat{\bb} \in \Nc(A^T)$ is obtained by the MATLAB function, e.g., {\sf null}, and $\delta$ is the noise level. The test starts from two initial vectors $\bx^{(0)}=\textsc{0}$ and $\by^{(0)}=\bb$. The performance is compared in terms of iteration number, relative solution error, and computing time in seconds, which are respectively abbreviated as IT, RSE, and Time(s). RSE is defined by RSE $= \bT{\bx^{(k)} - \bx^{\ast}}/\bT{\bx^{\ast}}$ for $k=0,1,2,\cdots$ and Time(s) is measured by the MATLAB built-in function {\sf tic-toc}.

\begin{figure}[!htbp]
\centering
    \subfigure[$n=100$, $m=30n$, $r=n/2$, $\kappa=n/10$]{
		\includegraphics[width=0.5\textwidth]{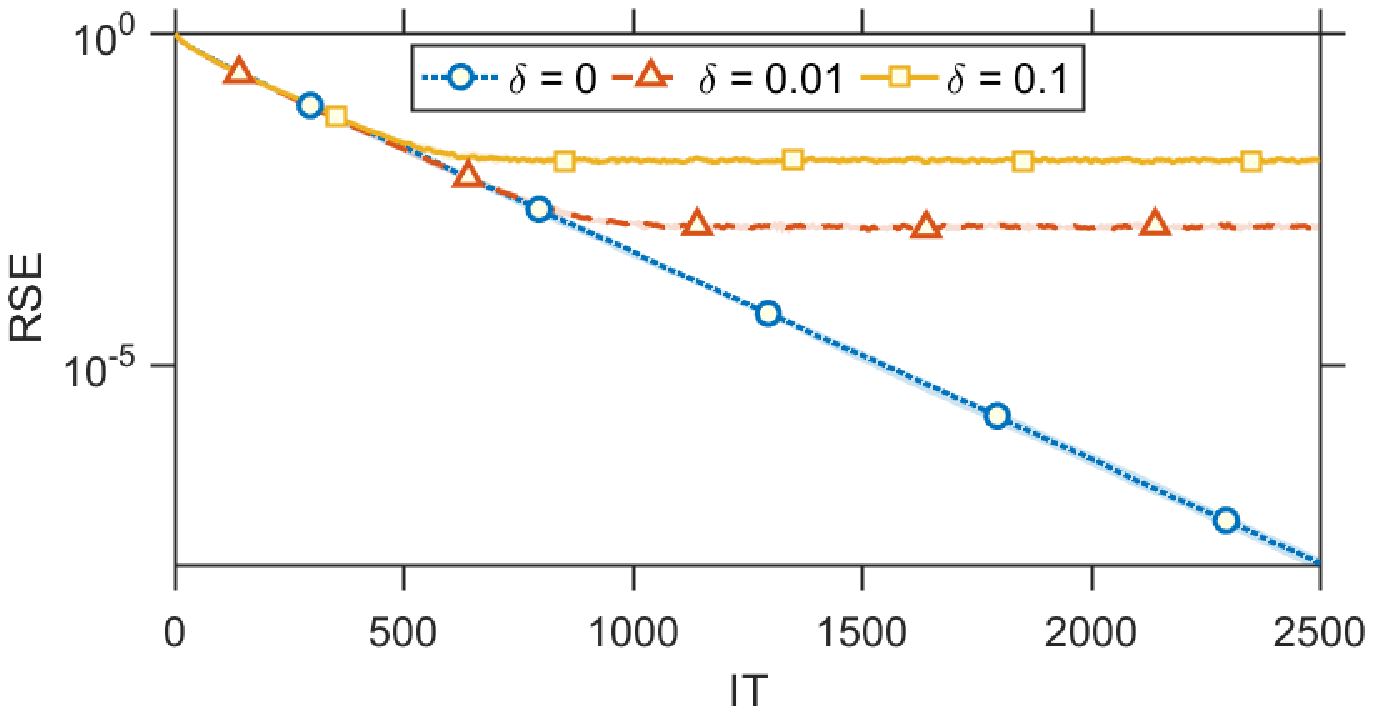}}
    \hspace{-0.7cm}
    \subfigure[$n=200$, $m=30n$, $r=n/2$, $\kappa=n/10$]{
		\includegraphics[width=0.5\textwidth]{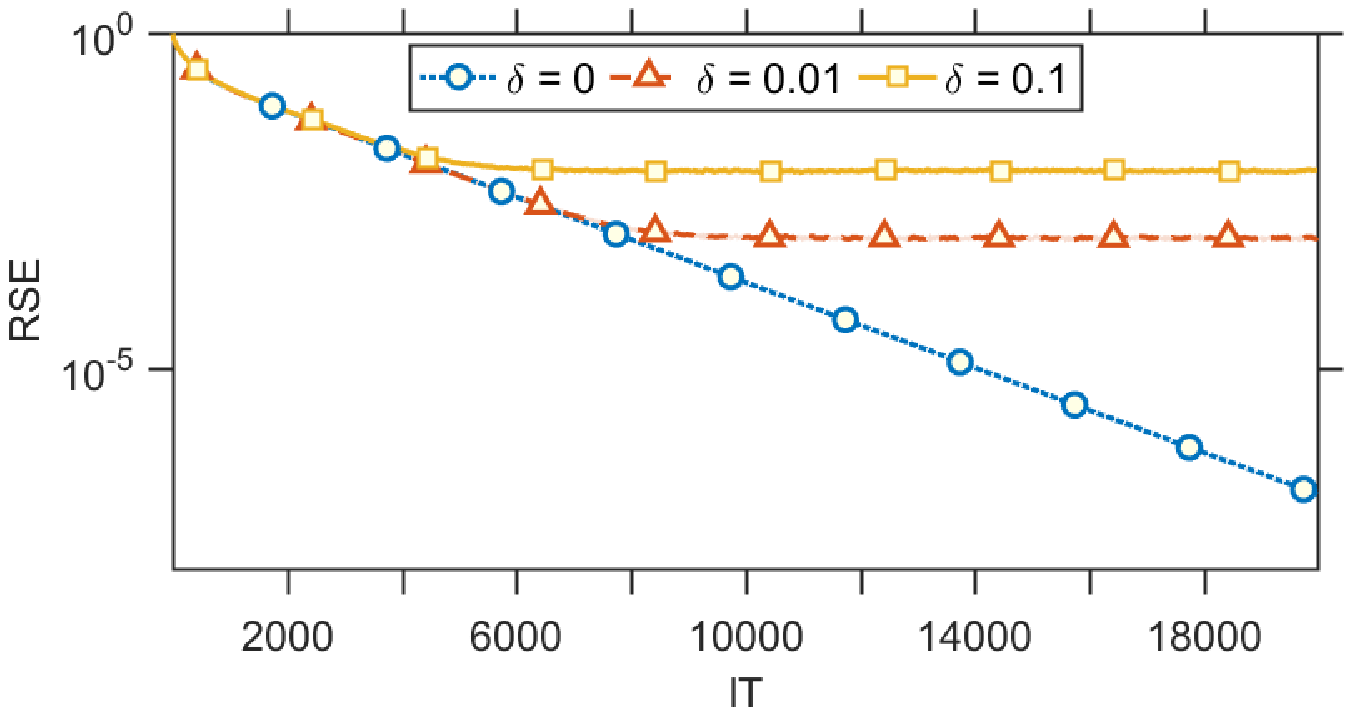}}
\caption{The convergence behaviors of RSE versus IT given by RMR for Example \ref{ERMR:example1} with $n=100$ (left) and $200$ (right) when the noise levels $\delta = 0$, $0.01$, and $0.1$.}
\label{fig:RMR:example1_SemiConverge}
\end{figure}

To begin with we depict the convergence behaviors of RSE versus IT, given by RMR, in Figure \ref{fig:RMR:example1_SemiConverge} with $\delta=0$, $0.01$, and $0.1$. In this figure, a heavy line represents median performance and the shaded region spans the minimum to the maximum value across all trials after repeatedly running this randomized iterative method $10$ times.  When $\delta=0$, the system \eqref{eq:Ax=b} is consistent. In this setting, we can see that RMR converges to the least-squares solutions successfully. When $\delta=0.01$ and $0.1$, the system \eqref{eq:Ax=b} is inconsistent. It is shown that the RMR iteration initially converges toward a good approximation of $\bx^{\ast}$ while continuing the iteration leads to corruption in the iteration vectors by noise and the approximation starts to stagnate and hovers around $\bx^{\ast}$. This implies that RMR has the semi-convergence behavior when solving the inconsistent linear system, and the number of iteration steps acts as a regularization parameter.

\begin{figure}[!htbp]
\centering
    \subfigure[$n=100$, $m=30n$, $r=n/2$, $\kappa=n/10$]{
		\includegraphics[width=0.5\textwidth]{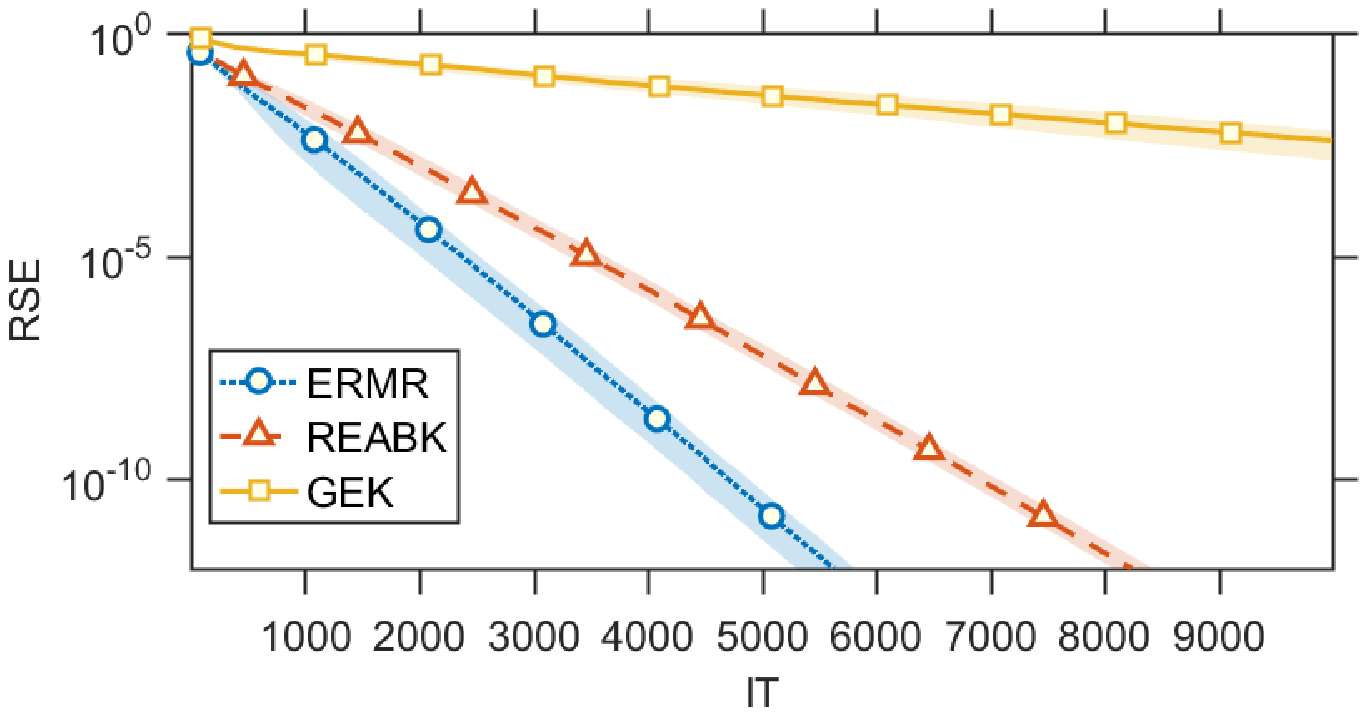}}
    \hspace{-0.7cm}
    \subfigure[$n=100$, $m=30n$, $r=n/2$, $\kappa=n/10$]{
		\includegraphics[width=0.5\textwidth]{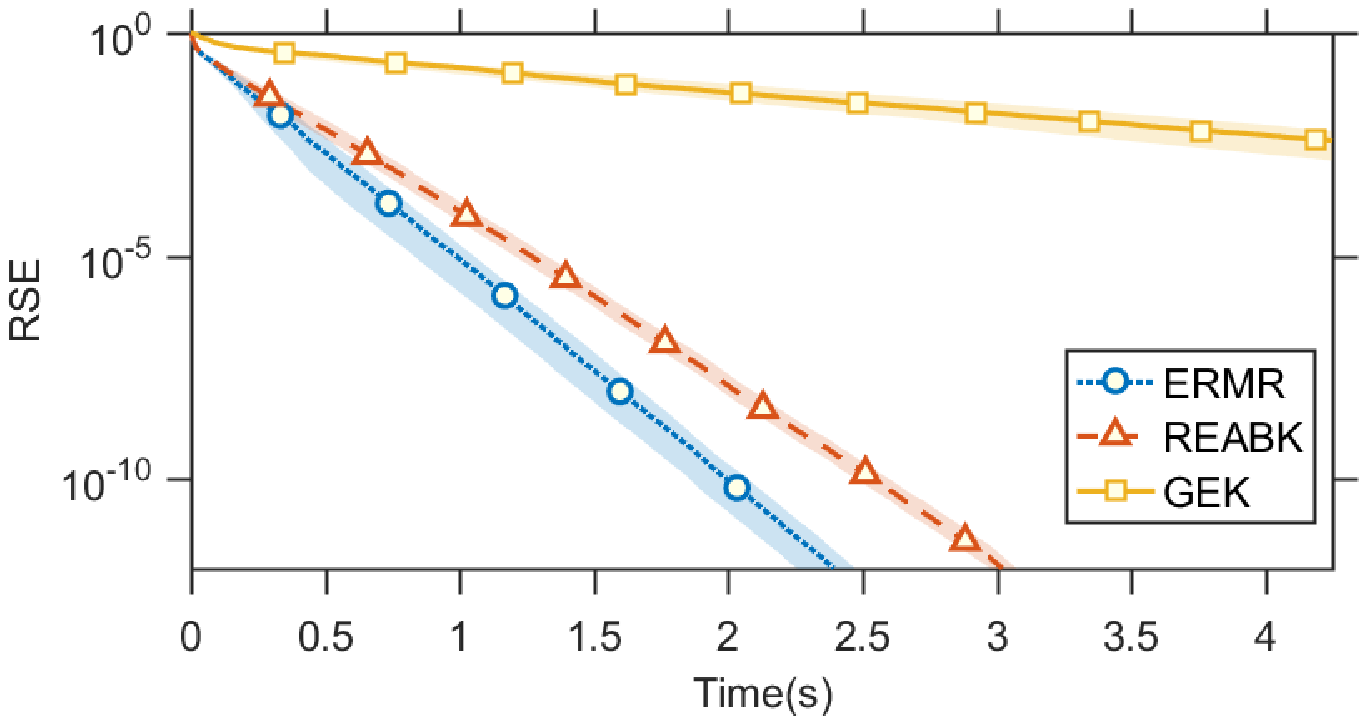}}
    \hspace{3cm}
    \subfigure[$n=200$, $m=30n$, $r=n/2$, $\kappa=n/10$]{
		\includegraphics[width=0.5\textwidth]{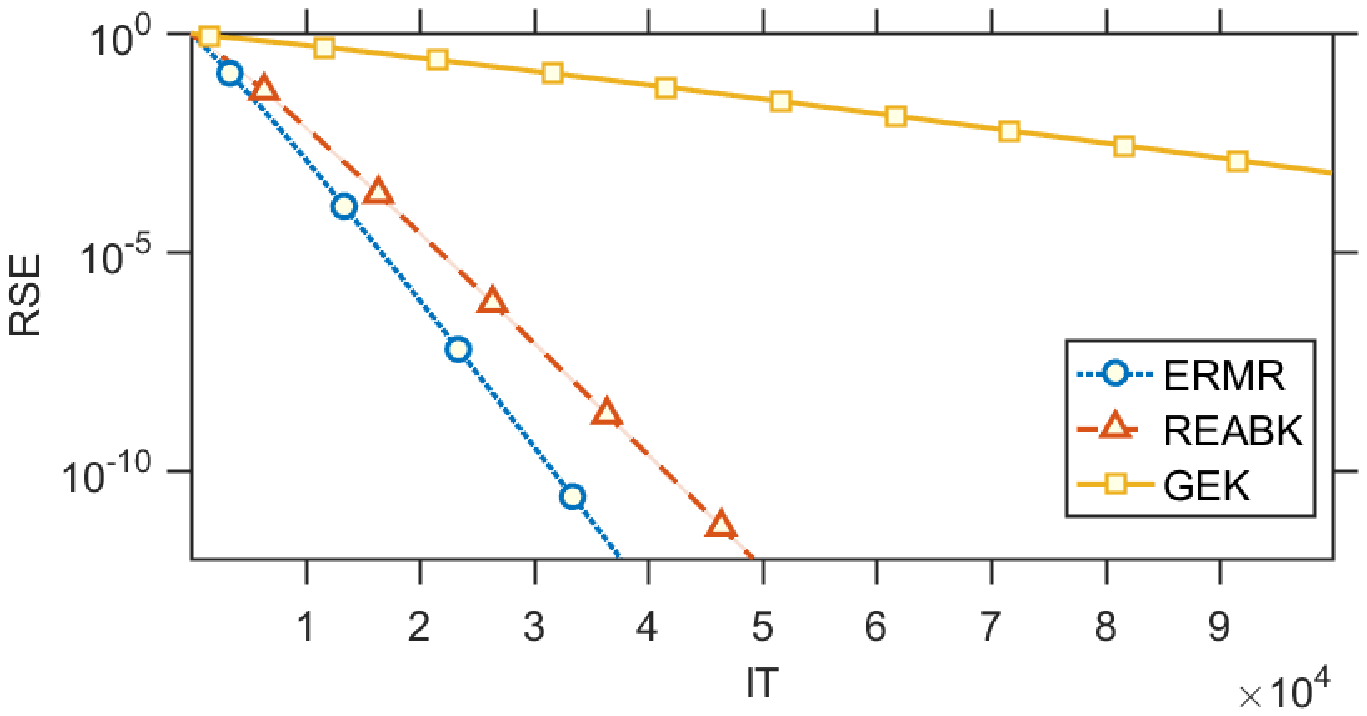}}
    \hspace{-0.7cm}
    \subfigure[$n=200$, $m=30n$, $r=n/2$, $\kappa=n/10$]{
		\includegraphics[width=0.5\textwidth]{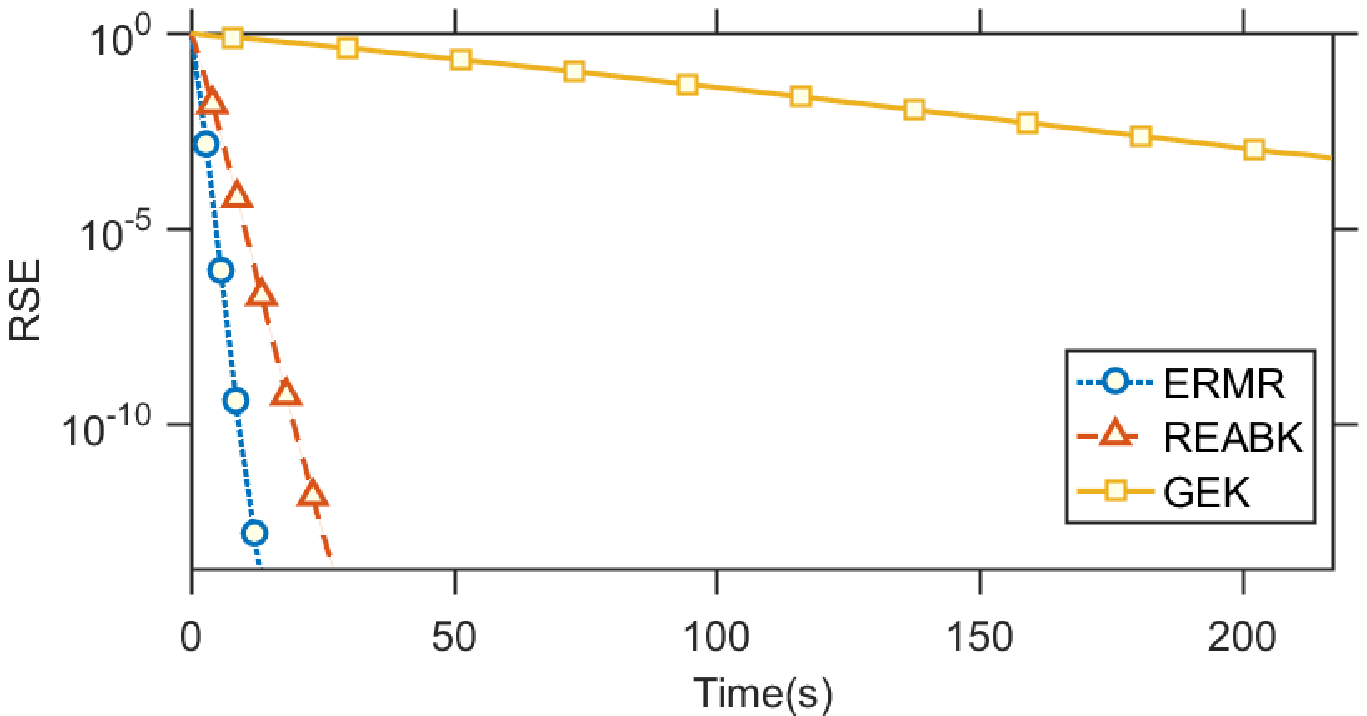}}
\caption{The convergence behaviors of RSE versus IT (left) and  Time(s) (right) given by ERMR, REABK, and GEK for Example \ref{ERMR:example1} with various $m$, $n$, $r$, and $\kappa$  when the noise level $\delta=0.1$.}
\label{fig:ERMR:example1}
\end{figure}

Besides, we track the convergence behaviors of RSE versus IT and  Time(s), given by ERMR, REABK, and GEK, and compare their performances for solving these same inconsistent linear systems. We take the noise level $\delta=0.1$ as an example. The results are shown in Figure \ref{fig:ERMR:example1}. We can observe the following phenomena. (I) The convergence curves of ERMR, REABK, and GEK break the semi-convergence horizons of RMR,  and  all the tested methods converge to the least-squares solution linearly. (II) The relative solution error of ERMR decay faster than that of REABK  and much faster  than that of GEK when the number of iteration steps and computing time increase.

\subsection{Image reconstruction} In the following experiments, we investigate some numerical tests taken from the field of tomography image reconstruction.

\vskip 0.5ex
\begin{example}\label{AIR:Exp1}
We solve the least-squares problem chosen from AIR Tools II toolbox for MATLAB (available from
\href{https://github.com/jakobsj/AIRToolsII}{https://github.com/jakobsj/AIRToolsII.}) \cite{18HJ}. We select two types of test problems: $2$-dimension seismic travel-time tomography and spherical Radon  transform tomography. They are represented by $[A,\bb, \bx^{\ast}] = {\sf seismictomo}(N, s, p)$ and ${\sf sphericaltomo}(N, \theta, p)$, which generate a coefficient matrix $A$, an exact solution $\bx^{\ast}$, and an exact right-hand side $\bb=A \bx^{\ast}$ by adjusting the input parameters, such as the size of the discrete domain ($N$),  the number of sources  ($s$), the projection angles in degrees ($\theta$), and the  number of rays ($p$).
\end{example}

\vskip 0.5ex
This example presents a sparse matrix $A$ with a full column-rank.  In particular, $s$ equals the length of $\theta$. The resulting matrix is of size $(sp)\times(N^2)$. In this setting, the inconsistent linear system is realized by setting the noisy right-hand side as $\bb = A\bx^{\ast} + \widehat{\bb}/\bT{\widehat{\bb}}$, where $\widehat{\bb}$ is a nonzero vector in the null space of $A^T$ generated by {\sf null}. From the illustration of the measurement geometries in the test problems \cite{18HJ}, we know that a loop displays the rows of $A$ reshaped into an $N\times N$ image. Then, we fix the block sizes of ERMR and REABK as $\tau=N$.

\begin{figure}[!htbp]
\centering
    \subfigure[Exact phantom]{
		\includegraphics[width=0.35\textwidth]{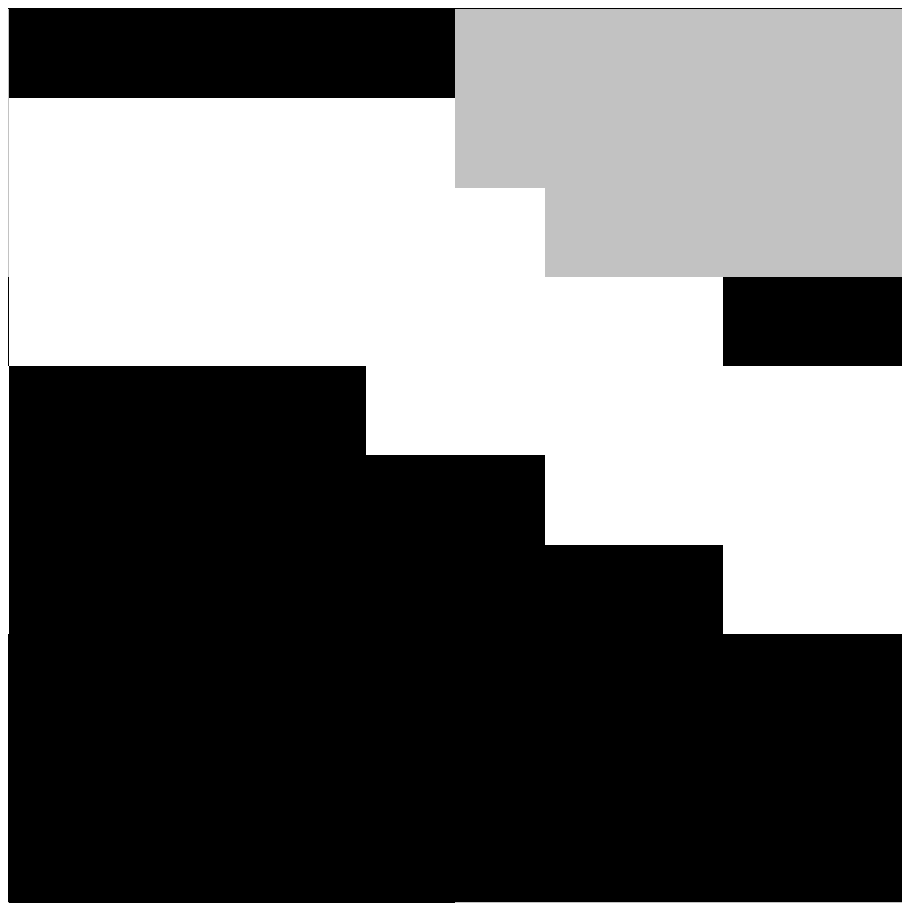}}
    \hspace{-1cm}
    \subfigure[Final phantom by REABK]{
		\includegraphics[width=0.35\textwidth]{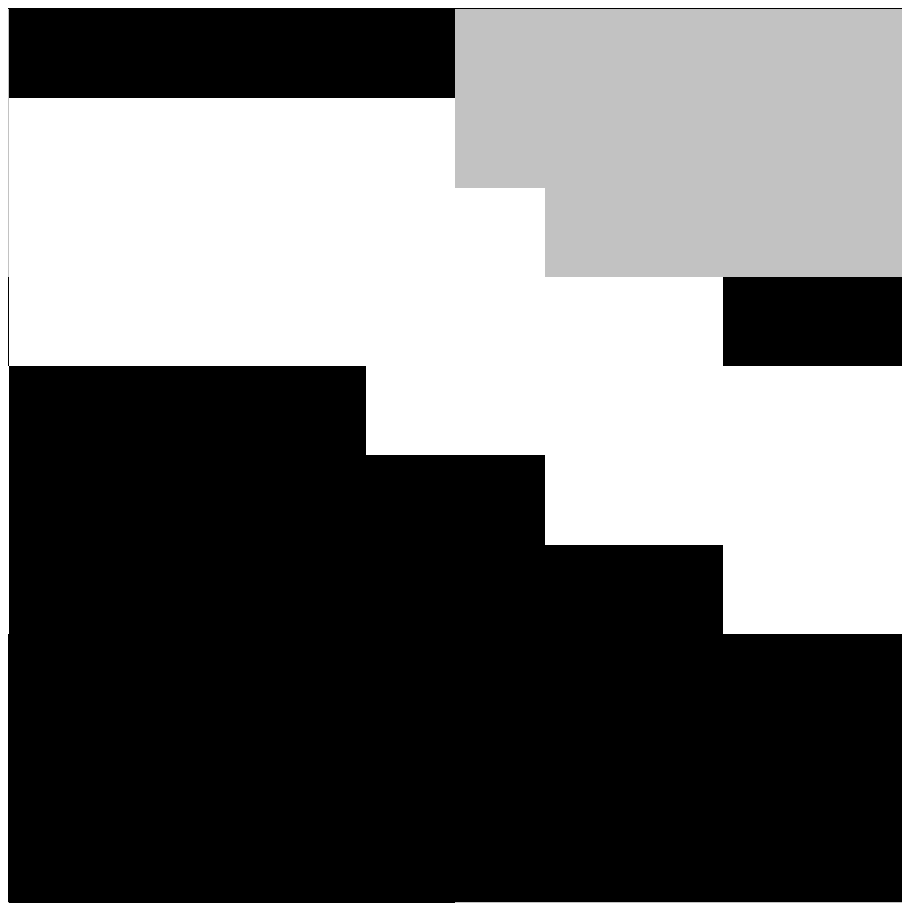}}
    \hspace{-1cm}
    \subfigure[Final phantom by ERMR]{
		\includegraphics[width=0.35\textwidth]{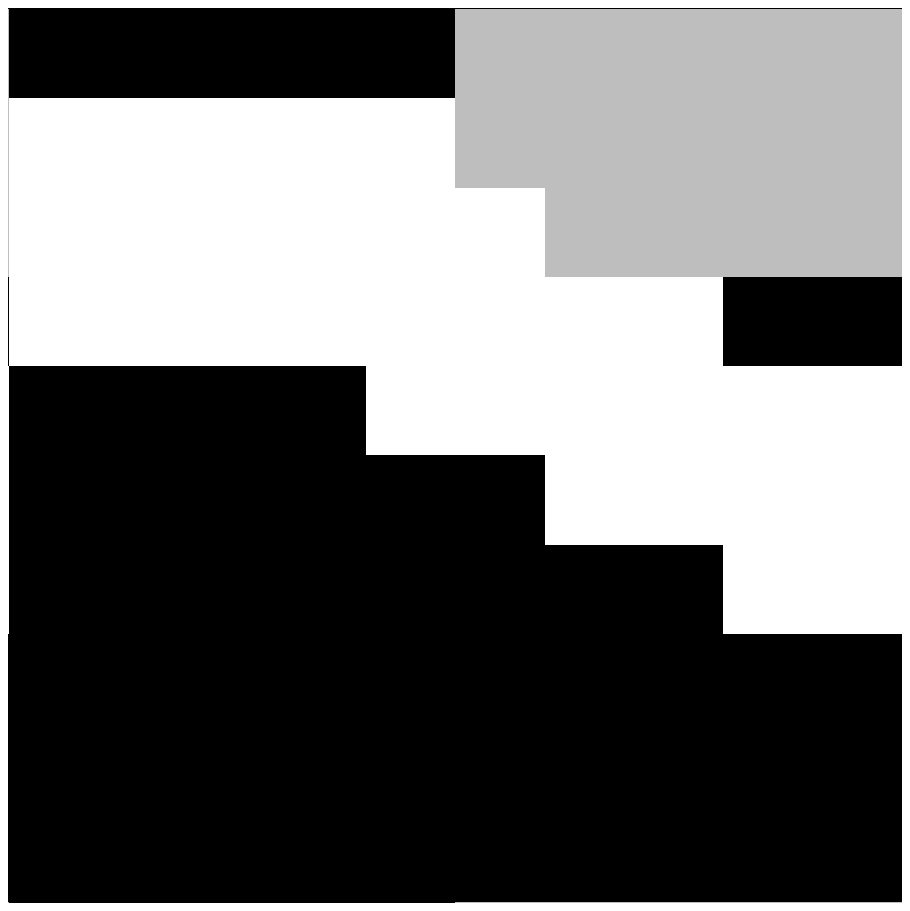}}
    \caption{Images of the exact phantom (a), and the approximate solutions obtained by REABK (b) and ERMR (c) for the seismic travel-time tomography.}
\label{fig:AIR:Exp1}
\end{figure}

In the seismic travel-time tomography test problem, we assign values to the input parameters as $N=10$, $s=180$, and $p=30$. The numerical results of RSE, IT, and Time(s), given by ERMR, REABK, and GEK, are listed in Table \ref{tab:AIR:Exp1}. It shows that ERMR and REABK successfully compute an approximate solution, but GEK fails due to the number of iteration steps exceeding $2\times 10^{6}$. For the convergent cases, the iteration counts and computing times of ERMR are  appreciably smaller than those of REABK. Hence, ERMR considerably outperforms REABK in terms of both iteration counts and computing times. In Figure \ref{fig:AIR:Exp1}, we give the $N \times N$ images of the exact tectonic phantom, and the approximate solutions obtained by ERMR and REABK. We see that the images constructed by them converge to the exact solution accurately.

\begin{table}[!ht]
 \normalsize
\caption{The numerical results of RSE, IT, and Time(s) obtained by ERMR, REABK, and GEK for the seismic travel-time tomography.}
\centering
\begin{tabular}{ cccc}
\cline{1-4}
&GEK&REABK&ERMR\\
\cline{1-4}
RSE      &$\#$             &$1.006\times 10^{-6}$&$9.169\times 10^{-7}$\\
IT       &$>2\times 10^{6}$&$1.398\times 10^{6} $&$2.928\times 10^{5} $\\
Time(s)  &$\#$             &441.3                &101.7\\
\cline{1-4}
\end{tabular}
\begin{tablenotes}
\footnotesize
\item[1.] {\it The item ' $>2\times 10^{6}$' represents that the number of iteration steps exceeds $2\times 10^{6}$. In this case, the corresponding RSE and Time(s) are expressed by $'\#'$.}
\end{tablenotes}
\label{tab:AIR:Exp1}
\end{table}

\begin{figure}[!htbp]
\centering
    \subfigure[Exact phantom]{
		\includegraphics[width=0.35\textwidth]{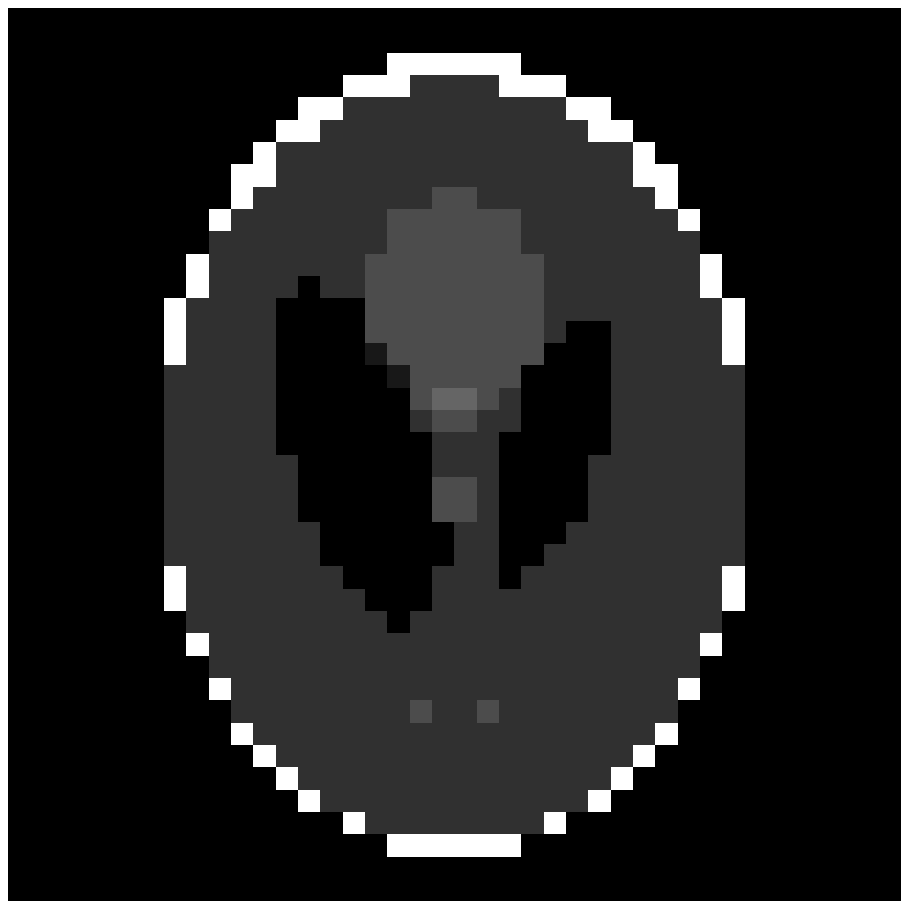}}
    \hspace{-1cm}
    \subfigure[Final phantom by REABK]{
		\includegraphics[width=0.35\textwidth]{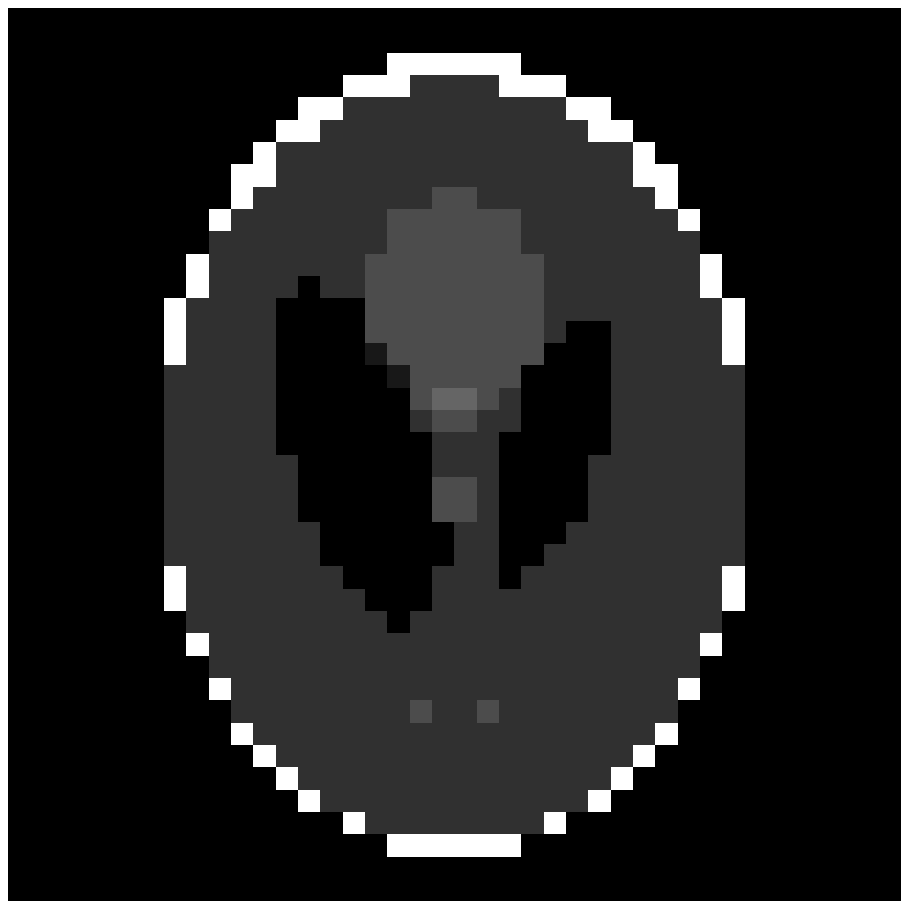}}
    \hspace{-1cm}
    \subfigure[Final phantom by ERMR]{
		\includegraphics[width=0.35\textwidth]{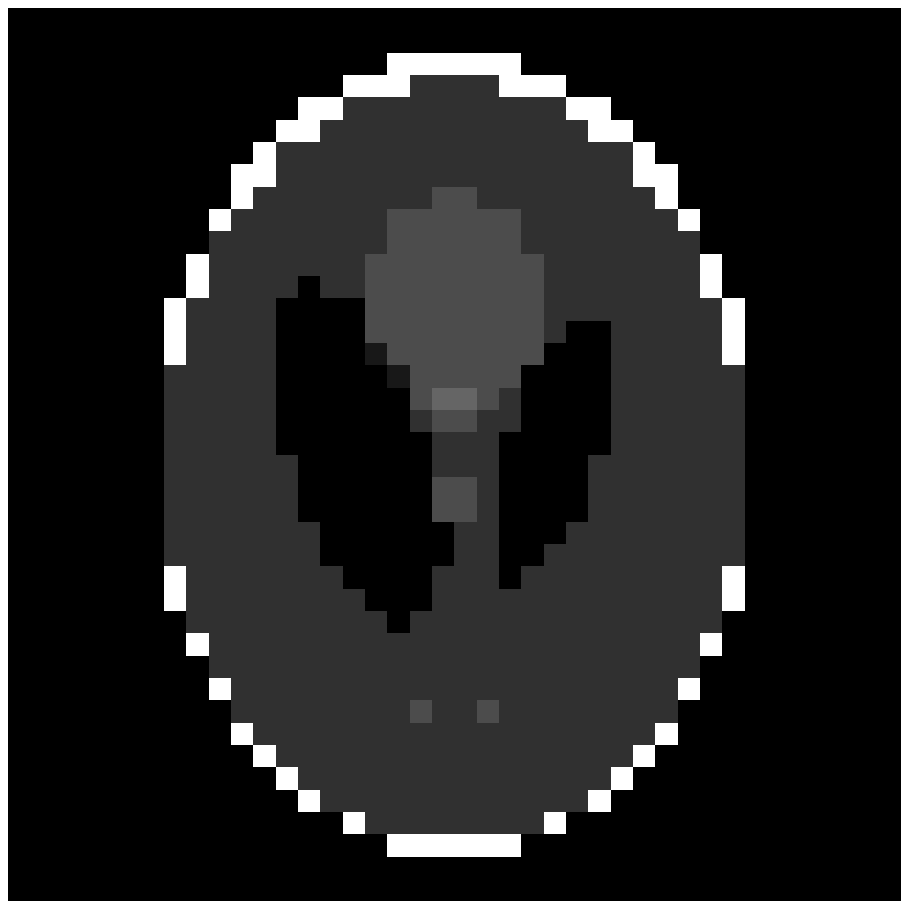}}
    \caption{Images of the exact phantom (a), and the approximate solutions obtained by REABK (b) and ERMR (c)  for the spherical Radon transform tomography.}
\label{fig:AIR:Exp2}
\end{figure}

In the spherical Radon transform tomography test problem,  the input parameters are assigned by $N=40$, $\theta=0:1:359$, and $p=20$. We list the numerical results of RSE, IT, and Time(s), given by ERMR, REABK, and GEK, in Table \ref{tab:AIR:Exp2}. We can also see that ERMR
 outperforms REABK and GEK in terms of both iteration counts and computing times. The exact tectonic phantom and the approximate solutions computed by ERMR and REABK are shown in Figure \ref{fig:AIR:Exp2}. We observe that these two methods successfully construct the exact image.

\begin{table}[!ht]
 \normalsize
\caption{The numerical results of RSE, IT, and Time(s) obtained by ERMR, REABK, and GEK for the spherical Radon transform tomography.}
\centering
\begin{tabular}{ cccc}
\cline{1-4}
&GEK&REABK&ERMR\\
\cline{1-4}
RSE      &$\#$             &$1.000\times 10^{-6}$&$9.996\times 10^{-7}$\\
IT       &$>2\times 10^{6}$&$4.857\times 10^{5} $&$2.542\times 10^{5} $\\
Time(s)  &$\#$             &2207.9                &989.7\\
\cline{1-4}
\end{tabular}
\label{tab:AIR:Exp2}
\end{table}

\subsection{Noisy data fitting} The subsequent numerical experiments consider the linear systems in CAGD research.

Rather than the traditional fitting methods directly solving a linear system, the geometric iterative method (GIM) generates a series of fitting curves by a fixed parameter \cite{18LMD}. Let us fit the ordered point set $ \{{\bm q}_i \in \Rc^3 : i\in[m] \}$.  Assume that $\left\{\mu_j(x):j\in [n]\right\}$ is a basis sequence and ${\bm p}_j^{(k)}$ is the $j$th control point at $k$th iterate. The $k$th GIM curve,
\begin{align}\label{eq:curve_C(x)}
  \mathcal{C}^{(k)}(x) =\sum_{j=1}^n \mu_j(x) {\bm p}_j^{(k)},~x \in [x_1,~x_m],
\end{align}
progressively approximates a target curve by updating the control points according to
\begin{align*}
 {\bm p}_j^{(k+1)} = {\bm p}_j^{(k)} +  {\bm \delta_j^{(k)}},
\end{align*}
where ${\bm \delta_j^{(k)}} \in \Rc^3$ is called the adjust vector and computed by ${\bm q}_i$ and ${\bm p}_j^{(k)}$. Let the $x$-, $y$-, and $z$-coordinates of the control point ${\bm p}_j^{(k)}$ (resp., the adjust vector ${\bm \delta}_j^{(k)}$) be respectively stored in the vectors ${\bm p}_x^{(k)}$, ${\bm p}_y^{(k)}$, and ${\bm p}_z^{(k)}$  (resp., ${\bm \delta}_x^{(k)}$, ${\bm \delta}_y^{(k)}$, and ${\bm \delta}_z^{(k)}$). From algebraic aspects, the GIM iterative processes,
\begin{align*}\renewcommand\arraystretch{1.5}
   {\bm p}_x^{(k+1)} = {\bm p}_x^{(k)} + {\bm \delta}_x^{(k)},~
   {\bm p}_y^{(k+1)} = {\bm p}_y^{(k)} + {\bm \delta}_y^{(k)},~{\rm and}~~
   {\bm p}_z^{(k+1)} = {\bm p}_z^{(k)} + {\bm \delta}_z^{(k)},
\end{align*}
are equal to solving three linear systems. Therefore, ERMR is suitable to deal with such problems.

The implementation detail of ERMR curve fitting is presented as follows. Let the data points be arranged into
${\bm q} = \left[{\bm q}_{1}~{\bm q}_{2}~\cdots ~{\bm q}_{m}\right]^T =\left[ {\bm q}_x ~ {\bm q}_y ~ {\bm q}_z\right] \in \Rc^{m\times 3}$. We input
the collocation matrix $A$,
two partitions of $[n]$,
initial vectors $\bm{p}_x^{(0)}\in \Rc^{n}$ (resp., $\bm{p}_y^{(0)}$, $\bm{p}_z^{(0)}$)
and $\bm{y}_x^{(0)}\in \Rc^{n}$ (resp., $\bm{y}_y^{(0)}$, $\bm{y}_z^{(0)}$),
the right-hand side $\bm{q}_x\in \Rc^{m}$ (resp., $\bm{q}_y$, $\bm{q}_z$),
and compute the next vector $\bm{p}_x^{(k)}$ (resp., $\bm{p}_y^{(k)}$, $\bm{p}_z^{(k)}$) using the ERMR update rule in Algorithm \ref{alg:ERMR}. Then, the approximate curve is formulated according to formula \eqref{eq:curve_C(x)}.

\vskip 0.5ex
\begin{example}\label{ERMR:example2}
Consider the least-squares problems in data fitting; see, e.g., the survey in \cite{18LMD}. We fit the data points $ \{{\bm q}_i: i\in[m] \}$ sampled from the granny knot curve, whose parametric equation is given as follows (available from \href{http://paulbourke.net/geometry/}{http://paulbourke.net/geometry/}).
$\bm{x} = -22\cos(t) - 128 \sin(t) - 44 \cos(3t) - 78 \sin(3t)$,
$\bm{y} = -10 \cos(2t) - 27 \sin(2t) + 38 \cos(4t) + 46 \sin(4t)$, and
$\bm{z} = 70 \cos(3t) - 40 \sin(3t)$ for $0 \leq t \leq 2\pi$.
As other researchers do, we first assign a parameter sequence $\bm{\nu}$ and a knot vector $\bm{\mu}$ of cubic B-spline basis, which is  simple and has a wide range of applications in CAGD; see, e.g., \cite{14DL}, and then obtain the collocation matrix using the MATLAB built-in function as $A={\sf spcol}(\bm{\mu}, 4, {\bm \nu})$.
\end{example}

\vskip .5ex
For more details on the formulations of generating $\bm{\nu}$ and $\bm{\mu}$, we respectively refer to equations (9.5) and (9.69) in the book by Piegl and Tiller \cite{97PT}. Though ERMR can be started with arbitrary initial control points, the works in \cite{18LCZ,18LMD} suggest that equation (23) in \cite{14DL} is an appropriate and effective selection strategy.  At $k$th iterate, the relative solution error is defined by RSE = $E_k/E_0$, where
$E_k = \bF{ {\bm p}^{(k)} -  {\bm p}^{\ast}}$ and ${\bm p}^{\ast} = A^{\dag} {\bm q}$ is the least-square solution.

\begin{figure}[!htbp]
\centering
    \subfigure[Initial data points]{
		\includegraphics[width=0.35\textwidth]{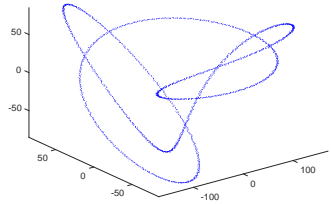}}
    \hspace{1cm}
    \subfigure[Final curve by GEK]{
		\includegraphics[width=0.35\textwidth]{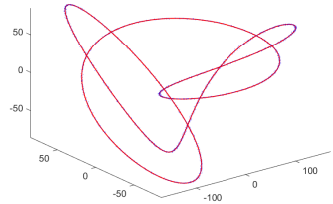}}
    \hspace{1cm}
    \subfigure[Final curve by REABK]{
		\includegraphics[width=0.35\textwidth]{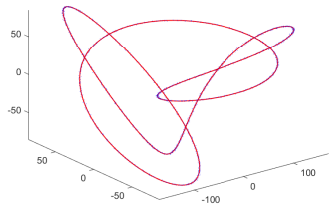}}
    \hspace{1cm}
    \subfigure[Final curve by ERMR]{
		\includegraphics[width=0.35\textwidth]{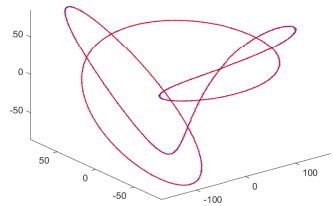}}
\caption{The initial data points (a) and the limit curves given by GEK(b), REABK(c), and ERMR(d) for Example \ref{ERMR:example2} when $n=200$ and $m=3000$.}
\label{fig:granny knot_curve}
\end{figure}

In the following, we discuss the capability of ERMR to fit the noisy data ${\bm q} = A{\bm p}^{\ast} + \widehat{{\bm q}}/\bF{\widehat{\bm q}}$, where $\widehat{\bm q}$ is a nonzero vector in the null space of $A^T$ generated by {\sf null}. The noisy initial data points are shown in Figure \ref{fig:granny knot_curve} (a) with $m=3000$ as a example. For solving this least-squares fitting problem,  we use $n$ control points to fit $m=15n$ data points and list the RSE, IT, and Time(s) given by ERMR, REABK, and GEK in Table \ref{tab:ERMR:example2+granny-knot}  when $n=200$, $400$, and $600$. We can see that ERMR outperforms REABK and GEK in terms of both iteration counts and computing times. In Figure \ref{fig:granny knot_curve} (b)-(d), we plot the limit curves constructed by ERMR, REABK, and GEK. It implies that these three methods achieve success in converging to the least-squares solution.

\begin{table}[!ht]
 \normalsize
\caption{The numerical results of RSE, IT, and Time(s) obtained by ERMR, REABK, and GEK for the granny knot curve fitting with $m=15n$.}
\centering
\begin{tabular}{ ccccc}
\cline{1-5}
&&GEK&REABK&ERMR\\
\cline{1-5}
$n=200$&RSE      &$1.000\times 10^{-6}$ &$9.689\times 10^{-7}$ &$9.780\times 10^{-7}$\\
&IT       &$1.821\times 10^{4}$ &$1.612\times 10^{4} $ &$1.129\times 10^{4} $\\
&Time(s)  &$21.554$              &$11.286$                 &$9.578$\\
\cline{1-5}
$n=400$&RSE      &$9.941\times 10^{-7}$ &$9.869\times 10^{-7}$ &$9.950\times 10^{-7}$\\
&IT       &$3.138\times 10^{4}$ &$2.817\times 10^{4} $ &$2.212\times 10^{4} $\\
&Time(s)  &$329.167$              &$35.485$                 &$26.400$\\
\cline{1-5}
$n=600$&RSE      &$1.000\times 10^{-6}$ &$9.675\times 10^{-7}$ &$9.950\times 10^{-7}$\\
&IT       &$5.175\times 10^{4}$ &$4.792\times 10^{4} $ &$3.312\times 10^{4} $\\
&Time(s)  &$1332.485$              &$77.170$                &$57.458$\\
\cline{1-5}
\end{tabular}
\label{tab:ERMR:example2+granny-knot}
\end{table}

\section{Conclusions}\label{sec:conclusions}
In this paper, the ERMR iterative method is presented for solving the inconsistent linear system.  ERMR uses two randomized multiple row iterations in each step and finds the Moore-Penrose inverse solution, i.e., the least-squares solution with minimum Euclidean-norm. It is shown that the proposed method has a linear convergence rate in the mean square.  We also provide the computational complexity analysis for ERMR. Some numerical examples, including the synthetic data and real-world  applications, are given to demonstrate the convergence behavior of ERMR when the coefficient matrix is dense or sparse, full rank or rank deficient. Numerical results show that ERMR provides significant computational advantages compared with the existing extended pseudoinverse-free block row iterative methods. It means that ERMR is a competitive row-type variant for solving the inconsistent linear system.

\section{Appendix} Before we present the proofs of Theorems \ref{LemmaRMR} and \ref{ERMR:thm+RMR+Ay=0} for completeness, we may need the following lemmas, whose proofs can be found in \cite[Lemmas 2.4 and 2.5]{20DSS}.

\vskip 0.5ex
\begin{lemma}\label{ERMR:LemmaSigmaMin}
  Let $A\in \Rc^{m\times n}$ be any nonzero real matrix. For every $\bu \in \Rc(A)$, it holds that
    $\bT{A^T \bu}\geq \sigma_{\min}(A)\bT{\bu}$, where $\sigma_{\min}(A)$ denotes the smallest nonzero singular value of $A$.
\end{lemma}

\vskip 0.5ex
\begin{lemma}\label{ERMR:LemmaSigmaMax}
  Let $A\in \Rc^{m\times n}$ be any nonzero real matrix. For every $\bu \in \Rc^{m}$, it holds that
    $\bT{AA^T \bu}\leq \sigma_{\max}(A)\bT{A^T\bu}$, where $\sigma_{\max}(A)$ denotes the largest nonzero singular value of $A$.
\end{lemma}

\vskip 1ex
\begin{proof} [of Theorem \ref{LemmaRMR}]
For $k=0,1,2,\cdots$, the iteration in Algorithm \ref{alg:RMR} satisfies that
  \begin{align*}
    \bx^{(k+1)}- \bx^{\ast}
    & = \bx^{(k)}- \bx^{\ast} + \frac{\bet_k^T(\bb -  A\bx^{(k)})}{\BT{A^T \bet_k}} A^T \bet_k\\
    & = \left( \bx^{(k)}- \bx^{\ast} \right) + \frac{\bet_k^T A (\bx^{\ast} -  \bx^{(k)})}{\BT{A^T \bet_k}} A^T \bet_k\\
    & = \left(I_n - \frac{A^T \bet_k \bet_k^T A }{\BT{A^T \bet_k}}\right) \left(\bx^{(k)}- \bx^{\ast}\right).
  \end{align*}
It is easy to show that
  \begin{equation*}
  \BT{\bx^{(k+1)}- \bx^{\ast}} = \BT{\bx^{(k)}- \bx^{\ast}} - \frac{|\bet_k^T A \left(\bx^{(k)}- \bx^{\ast}\right)|^2}{\BT{A^T \bet_k}}.
  \end{equation*}
The expected squared error follows that
\begin{equation*}
\bE\left[  \BT{\bx^{(k+1)}- \bx^{\ast}} \right] =
\bE\left[  \BT{\bx^{(k)}- \bx^{\ast}}   \right] -
\bE\left[  \frac{|\bet_k^T A \left(\bx^{(k)}- \bx^{\ast}\right)|^2}{\BT{A^T \bet_k}} \right].
\end{equation*}
According to the following two properties that
  \begin{align*}
    \bet_k^T A \left(\bx^{(k)}- \bx^{\ast}\right)
    & = \sum_{i\in \mathcal{I}_i} (b_i - A_{i,:} \bx^{(k)}) \bmu_{i}^T A \left( \bx^{(k)}- \bx^{\ast} \right)\\
    & = \sum_{i\in \mathcal{I}_i} (b_i - A_{i,:} \bx^{(k)}) \bmu_{i}^T (A \bx^{(k)} - \bb)\\
    & = -\BT{ \bb_{\mathcal{I}_i} -   A_{\mathcal{I}_i, :} \bx^{(k)}}
  \end{align*}
and
  \begin{align*}
    \BT{A^T \bet_k}
     & = \BT{ \sum_{i\in \mathcal{I}_i} (b_i - A_{i,:} \bx^{(k)}) A^T \bmu_{i} }\\
     & = \BT{ A_{\mathcal{I}_i, :} (\bb_{\mathcal{I}_i} -   A_{\mathcal{I}_i, :} \bx^{(k)}) }\\
     & \leq \sigma_{\max}^2\left( A_{\mathcal{I}_i, :} \right) \BT{ \bb_{\mathcal{I}_i} -   A_{\mathcal{I}_i, :} \bx^{(k)}},
  \end{align*}
where the last inequality is from Lemma \ref{ERMR:LemmaSigmaMax}, it yields that
  \begin{align*}
   \frac{|\bet_k^T A \left(\bx^{(k)}- \bx^{\ast}\right)|^2}{\BT{A^T \bet_k}}
   & \geq \frac{1}{\sigma_{\max}^2\left( A_{\mathcal{I}_i, :} \right)} \BT{ \bb_{\mathcal{I}_i} -   A_{\mathcal{I}_i, :} \bx^{(k)}}\\
   & = \frac{ \BF{ A_{\mathcal{I}_i, :} } }{\sigma_{\max}^2\left( A_{\mathcal{I}_i, :}\right)}
   \frac{ \BT{ \bb_{\mathcal{I}_i}- A_{\mathcal{I}_i, :} \bx^{(k)}} }{ \BF{ A_{\mathcal{I}_i, :} } } \\
   & \geq \frac{1}{\beta_{\max}^{I}} \frac{ \BT{ \bb_{\mathcal{I}_i}- A_{\mathcal{I}_i, :} \bx^{(k)}} }{ \BF{ A_{\mathcal{I}_i, :} } }.
  \end{align*}
Let $\bE_k\left[\cdot\right]$ denote the conditional expectation conditioned on the first $k$ iterations in RMR. It follows that
  \begin{align*}
    \bE_k \left[ \frac{|\bet_k^T A \left(\bx^{(k)}- \bx^{\ast}\right)|^2}{\BT{A^T \bet_k}} \right]
    & \geq \frac{1}{\beta_{\max}^{I}}
           \sum_{i=1}^s \frac{\BF{A_{\mathcal{I}_i, :}}}{\BF{A}}
           \frac{ \BT{ \bb_{\mathcal{I}_i}- A_{\mathcal{I}_i, :} \bx^{(k)}} }{ \BF{ A_{\mathcal{I}_i, :} } } \\
    & = \frac{\BT{\bb - A\bx^{(k)}}}{\beta_{\max}^{I} \BF{A}}.
  \end{align*}
Since $\bx^{(0)}$, $\bx^{\ast} \in \Rc(A^T)$, it is easy to show $\bx^{(k)} - \bx^{\ast}\in \Rc(A^T)$ by induction. Lemma \ref{ERMR:LemmaSigmaMin} indicates that
  \begin{align*}
    \bE_k\left[  \frac{|\bet_k^T A \left(\bx^{(k)}- \bx^{\ast}\right)|^2}{\BT{A^T \bet_k}} \right]
     \geq \frac{ \sigma_{\min}^2(A)\BT{ \bx^{(k)} - \bx^{\ast} }}{\beta_{\max}^{I} \BF{A}}.
  \end{align*}
Taking expectation again gives
  \begin{equation*}
    \bE\left[ \frac{|\bet_k^T A \left(\bx^{(k)}- \bx^{\ast}\right)|^2}{\BT{A^T \bet_k}} \right]
    \geq  \frac{1}{\beta_{\max}^{I}} \frac{\sigma_{\min}^2(A)}{\BF{A}}
    \bE\left[\BT{ \bx^{(k)} - \bx^{\ast} } \right].
  \end{equation*}
Then, the recurrence yields the result from the expected squared error equality. \hfill
\end{proof}

\vskip .5ex

\begin{proof} [of Theorem \ref{ERMR:thm+RMR+Ay=0}]
This proof is similar to that of Theorem \ref{LemmaRMR} and yet the involved technicalities are somewhat different. For completeness and simplicity, we divide this proof into the following two parts.

First, the expected squared error follows that
\begin{equation*}
\bE\left[  \BT{\by^{(k+1)}-\bb_{N}} \right] =
\bE\left[  \BT{\by^{(k)}-\bb_{N}}   \right] -
\bE\left[  \frac{|\bze_k^T A^T\left(\by^{(k)}-\bb_{N}\right)|^2}{\BT{A \bze_k}}\right].
\end{equation*}

Second, based on several properties of $\bze_k$, including
\begin{align*}
    \bze_k^T A^T\left(\bb_{N} - \by^{(k)}\right) = \BT{A^T_{:,\mathcal{J}_j} \by^{(k)}}
~~ {\rm and} ~~
    \BT{A \bze_k} \leq \sigma_{\max}^2\left( A_{:,\mathcal{J}_j} \right)\BT{A^T_{:,\mathcal{J}_j}\by^{(k)}},
\end{align*}
it yields that
\begin{align*}
   \frac{|\bze_k^T A^T\left(\by^{(k)}-\bb_{N}\right)|^2}{\BT{A \bze_k}}
   \geq \frac{1}{\beta_{\max}^{J}} \frac{\BT{A^T_{:,\mathcal{J}_j}\by^{(k)}}}{\BF{A^T_{:,\mathcal{J}_j}}}.
\end{align*}
As a result, we have
  \begin{equation*}
    \bE\left[\frac{|\bze_k^T A^T\left(\by^{(k)}-\bb_{N}\right)|^2}{\BT{A \bze_k}}\right]
    \geq  \frac{1}{\beta_{\max}^{J}} \frac{\sigma_{\min}^2(A)}{\BF{A}}
    \bE\left[\BT{\by^{(k)}-\bb_{N} } \right].
  \end{equation*}
Then, the recurrence yields the result from the expected squared error equality. \hfill
\end{proof}

\section*{Acknowledgment}
This work is supported by the National Natural Science Foundation of China under grants 12101225 and 12201651.

\end{document}